\documentclass[10pt]{amsart}
\usepackage{amsfonts}
\usepackage{color}
\usepackage[square,compress,comma, numbers]{natbib}
\usepackage[colorlinks=true, citecolor=blue, linkcolor=blue]{hyperref}
\usepackage{amssymb}
\usepackage{amsmath}

\definecolor{c20}{rgb}{0.,0.7,0.}
\definecolor{c30}{rgb}{0.,0.,1.}
\definecolor{c40}{rgb}{1,0.1,0.7}
\definecolor{c50}{rgb}{1,0,0}
\definecolor{c60}{rgb}{1,0.9,0.1}

\allowdisplaybreaks[4]

\def\cEc#1{\textcolor{c30}{#1}}

\def\cEc#1{#1}

\newcommand{\kb}[1]{\boldsymbol{#1}}
\newcommand{\vk}[1]{\kb{#1}}

\newcommand{\abs}[1]{\left\lvert #1 \right\rvert}

\newcommand{\E}[1]{\mathbb{E}\left\{#1\right\}}

\newcommand{\pk}[1]{\mathbb{P} \left \{#1 \right \} }

\newcommand{\R}{\mathbb{R}}

\newcommand{\N}{\mathbb{N}}

\newcommand{\inn}{\in \N}
\newcommand{\ldot}{,\ldots,}

\newcommand{\BQN}{\begin{eqnarray}}
\newcommand{\EQN}{\end{eqnarray}}
\newcommand{\BQNY}{\begin{eqnarray*}}
\newcommand{\EQNY}{\end{eqnarray*}}

\newcommand{\BS}{\begin{sat}}
\newcommand{\ES}{\end{sat}}
\newcommand{\BT}{\begin{theo}}
\newcommand{\ET}{\end{theo}}
\newcommand{\BL}{\begin{lem}}
\newcommand{\EL}{\end{lem}}
\newcommand{\BK}{\begin{korr}}
\newcommand{\EK}{\end{korr}}

\newcommand{\BD}{\begin{de}}
\newcommand{\ED}{\end{de}}

\newcommand{\BIT}{\begin{itemize}}
\newcommand{\EIT}{\end{itemize}}
\newcommand{\BDI}{\begin{description}}
\newcommand{\EDI}{\end{description}}

\newcommand{\BRM}{\begin{remarks}}
\newcommand{\ERM}{\end{remarks}}

\newcommand{\BEL}{\begin{lem}}
\newcommand{\EEL}{\end{lem}}

\newtheorem{theo}{Theorem}[section]
\newtheorem{sat}[theo]{Proposition}
\newtheorem{de}[theo]{Definition}
\newtheorem{lem}[theo]{Lemma}

\newtheorem{example}[theo]{Example}
\newtheorem{korr}[theo]{Corollary}
\newtheorem{remark}[theo]{Remark}
\newtheorem{remarks}[theo]{Remarks}
\newtheorem{prop}[theo]{Proposition}

\newcommand{\nelem}[1]{{Lemma \ref{#1}}}

\newcommand{\netheo}[1]{{Theorem \ref{#1}}}

\newcommand{\prooftheo}[1]{ \textsc{\bf Proof of Theorem} \ref{#1}:}
\newcommand{\proofprop}[1]{\textsc{\bf Proof of Proposition} \ref{#1}:}
\newcommand{\prooflem}[1]{\textsc{\bf Proof of Lemma} \ref{#1}:}

\newcommand{\proofExamp}[1]{\textsc{\bf Proof of Eaxmple} \ref{#1}:}
\newcommand{\COM}[1]{}

\newcommand{\QED}{\hfill $\Box$}

\def\rw{\rightarrow}

\def\IF{\infty}

\def\LT{\left}
\def\RT{\right}

\def\rw{\rightarrow}

\def\X{\vk X}

\def\vn{\varepsilon}
\def\Var{\text{Var}}

\def\II{\mathbb{I}}

\def\Bu+#1{\mathcal{B}^{\varepsilon+}_{u}(#1)}

\begin{document}

\title{Extremes of $L^p$-norm of Vector-valued Gaussian processes with Trend}
\author{Long Bai}
\address{Long Bai,
Department of Actuarial Science, %\\Faculty of Business and Economics\\
University of Lausanne\\
UNIL-Dorigny, 1015 Lausanne, Switzerland
}
\email{Long.Bai@unil.ch}
\COM{

\author{Krzysztof D\c{e}bicki}
\address{Krzysztof D\c{e}bicki, Mathematical Institute, University of Wroc\l aw, pl. Grunwaldzki 2/4, 50-384 Wroc\l aw, Poland}
\email{Krzysztof.Debicki@math.uni.wroc.pl}

\author{Enkelejd  Hashorva}
\address{Enkelejd Hashorva, Department of Actuarial Science, %\\Faculty of Business and Economics\\
University of Lausanne,\\
UNIL-Dorigny, 1015 Lausanne, Switzerland
}
\email{Enkelejd.Hashorva@unil.ch}

\author{Lanpeng Ji}
\address{Lanpeng Ji, Department of Actuarial Science, %\\Faculty of Business and Economics\\
University of Lausanne\\
UNIL-Dorigny, 1015 Lausanne, Switzerland
}
\email{Lanpeng.Ji@unil.ch}}

\bigskip

%\date{\today}
 \maketitle

{\bf Abstract:}
Let $\vk{X}(t)=(X_1(t),\ldots,X_n(t))$ be a Gaussian vector process and $g(t)$ be a continuous function. The asymptotics of distribution of $\LT\|\vk{X}(t)\RT\|_p$, the $L^p$ norm for Gaussian finite-dimensional vector, have been investigated in numerous literatures. In this contribution we are concerned with the exact tail asymptotics of $\LT\|\vk{X}(t)\RT\|^c_p,\ c>0, $ with trend $g(t)$ over $[0,T]$. Both scenarios that $\vk{X}(t)$ is locally stationary and  non-stationary are considered.
Important examples include $\sum_{i=1}^n \abs{X_i(t)}+g(t)$ and chi-square processes with trend, i.e., $\sum_{i=1}^n X_i^2(t)+g(t)$. These results are of interest in applications in engineering, insurance and statistics, etc.



\def\Xi{\chi}
\section{Introduction}
In engineering sciences, extreme values of non-linear functions of multivariate Gaussian processes are of interest in dealing with the safety of structures, see \cite{EXCRA1980} and the references therein. Probabilistic structural analysis to answer the question is: what is the probability that a certain mechanical (or other) structure will survive when it is subject to a random load.
The load is then usually defined by some $n$-dimensional vector process $\vk{Y}(t)=(Y_1(t),\ldots,Y_n(t)),\ n\geq 1,\ t\in[0,T]$, and one seeks the probability that $\vk{Y}$ exceeds some more or less well-defined safe region, which is specific for the structure as
\BQN\label{AA}
\pk{\vk{Y}(t)\notin \vk{S}_u(t), \text{for some}\ t\in [0,T]},
\EQN
where the time-dependent safety region $\vk{S}_u(t)$ is defined by
\BQNY
\vk{S}_u(t)=\left\{(x_1,\cdots, x_n)\in\R^n:\|\vk{x}\|_p\leq h(t,u)\right\}
\EQNY
with $h(t,u), \ t,u\geq 0$ some continuous function and  $||\cdot||_p$, $p\in[1,\IF]$ the $L^p$ norm, i.e.,
$$||\vk{x}||_p=\LT\{
\begin{array}{ll}
\LT(\sum_{i=1}^n|x_{i}|^p\RT)^{1/p},&\ p\in[1,\IF),\\
\max(|x_{1}|,\ldots,|x_{n}|),&\ p=\IF,
\end{array}
\RT.$$
in the space $L^p_n=\{\vk{x}=(x_1,\ldots,x_n):||\vk{x}||_p<\IF\}$.\\
Assume that  $\vk{X}(t)=(X_1(t),\ldots,X_n(t))$ where $X_i's$ are independent copies of $X(t)$ a centered Gaussian process which has continuous trajectories,  variance function $\sigma^2(\cdot)$ and correlation function $r(\cdot,\cdot)$
and
 \BQN\label{vd}
 \vk{d}=(d_1,\ldots,d_n), \ \ 1=d_1=\cdots=d_m>d_{m+1}\geq d_{m+2}\geq\cdots\geq d_{n}>0, \quad 1\leq m\leq n.
 \EQN
In the framework of \eqref{AA}, set $\vk{Y}(t)=\vk{d}*\vk{X}(t):=(d_1X_1(t),\cdots,d_nX_n(t))$, then we can rewrite \eqref{AA} as
\BQNY
\pk{\exists_{t\in[0,T]}Z(t)>h(t,u)}
\EQNY
where
\BQN\label{ZZ}
Z(t):=Z_{p}(t):=\LT\|\vk{X}(t)*\vk{d}\RT\|_p,
\EQN
and hereafter, we call $Z_p(t)$ the \emph{$L^p$ norm process}.\\
When $p=2$, for a positive constant $c$, as in the convention $Z^c_2(t)=\LT(Z_2(t)\RT)^c$ is called the \emph{chi process} when $c=1$ and the \emph{chi-square process} when $c=2$.\\
Further, as the Gaussian processes, we can introduce the \emph{stationary, locally-stationary, and non-stationary $L^p$ norm processes} according to the stationary, locally-stationary, and non-stationary properties of $X(t)$, respectively.\\
The investigate of
$$\pk{\exists_{t\in[0,T]}Z_2(t)>u}=\pk{\sup_{t\in[0,T]}Z_2(t)>u},\quad \text{as} \ u\rw\IF$$
is initiated by the studies of high excursions of envelope of a Gaussian process, see e.g., \cite{BN1969} and generalized in \cite{Lindgren1980a,Lindgren1980b,Lindgren1989}. When $X(t)$ is stationary with $\sigma(t)\equiv 1$ and
\BQNY
r(s,t)=1-a|t-s|^{\alpha} +o(|t-s|^{\alpha}),\ \ \abs{t-s}\rw 0,\ \alpha\in(0,2],
\EQNY
\cite{Albin1990,Albin1992} develop the Berman's approach in \cite{Berman82} to obtain an asymptotic behavior of large deviation probabilities of the stationary chi-square processes. \\
Further, if there exists unique $t_0\in[0,T]$ satisfies $\sigma(t_0)=\sup_{t\in[0,T]}\sigma(t)$ and
\BQNY
\sigma(t)=1-b(t_0)\abs{t-t_0}^2+o(\abs{t-t_0}^2),\quad
r(s,t)=1-a(t_0)\abs{t-s}^2+o(\abs{t-s}^2), \quad s,t\rw t_0,
\EQNY
where $b(t_0)$ and $a(t_0)$ are positive constants related to $t_0$,
the tail asymptotic behavior of the non-stationary $Z^2_2(t)$ and $Z_p(t),p\in(1,2)\cup(2,\IF)$ are investigated in  \cite{Pitchi1994} and \cite{FatalovA1993}, respectively, under the application of the so-called "double-sum method" in \cite{Pit96}.

Some recent contributions are focused on more general scenarios of chi process and chi-square process with $h(t,u)=u-g(t)$, i.e.,
\BQNY
\pk{\exists_{t\in[0,T]}Z_2^c(t)>h(t,u)}=\pk{\sup_{t\in[0,T]}\LT(Z_2^c(t)+g(t)\RT)>u},\ c=1,2,
\EQNY
where the continuous function $g(t)$ is generally considered as a trend or a drift.\\
When $X_i, i=1\ldot n$ are \emph{non-stationary Gaussian processes}, $Z_2(t)+g(t)$, the non-stationary chi processes with trend,  and $Z^2_2(t)-wt^\beta, w,\beta>0$, the non-stationary chi-square processes with trend, are studied in \cite{EnkelejdJi2014Chi} and \cite{PL2015}, respectively.\\
When $X_i, i=1\ldot n$ are \emph{locally-stationary Gaussian processes}, \cite{LJ2017} obtains the extreme of the supremum of $Z_2^2(t)$ with trend, see, e.g., \cite{Berman92,Hus90} for more details about locally stationary Gaussian processes.

Considering both the locally stationary and non-stationary $L^p$ norm processes, the contribution of this paper concerns an exact asymptotic behavior of large deviation probabilities for $Z_p^c(t)+g(t)$ with $p\in[1,\IF]$, constant $c\in(0,\IF)$ and $g(t),\ t\in[0,T]$ a continuous function, which contains the aforementioned results.\\
Organisation of the rest of the paper: In Section 2, the notation and some preliminaries are given. Our main results are displayed in Section 3. Following in Section 4 are two applications related to insurance and statistics.  Finally, we
present the proofs in Section 5 and several lemmas in Section 6.

\COM{For $u>0$, consider a vector-valued Gaussian  process
$\{{\bf{X}}_u(t),t\ge 0\}$, where ${\bf{X}}_u(t)=(X_{u,1}(t) \ldot X_{u,n}(t))$
with $\{X_{u,i}(t),t\ge 0\}, 1\le i\le n$, $n\inn$, being independent centered  Gaussian processes
with almost surely (a.s.) continuous sample paths.
In this paper we focus on the asymptotic behavior of the
probability that ${\bf{X}}_u$ enters the upper  orthant
$\{(x_1 \ldot x_n): x_i>M(u), i\in \{1 \ldot n\}\}$ over a fixed time interval $[0,T]$ with $\lim_{u\rw\IF}\frac{M(u)}{u}=c_0\in(0,\IF)$,  i.e.,
\BQN \label{aim}
\pk{\exists_{t \in [0,T]}
 \forall_{i=1 \ldot n} X_{u,i}(t)>M(u)},\ \ u\to\infty ,
\EQN
which is equivalent to
\BQN\label{aim2}
 p_{T,u}:=\pk{\sup_{t\in [0,T]} \min_{1 \le i \le n} X_{u,i}(t) > M(u)},\ \ u\to\infty ,
\EQN
implying that one can view (\ref{aim}) as the probability of extremal
behavior of the  process  $\LT\{\min_{1 \le i \le n} X_{u,i}(t), t\ge 0\RT\}$.}

\section{Notation and preliminaries}
First we introduce some notation, starting with the well-known Pickands constant $\mathcal{H}_{\alpha}$ defined by
$$\mathcal{H}_{\alpha}=\lim_{S_1\rightarrow\IF}\frac{1}{S_1}\mathcal{H}_{\alpha}[0,S_1],\quad \text{with }
\mathcal{H}_{\alpha}[-S_1,S_2]=\E{\sup_{t\in[-S_1,S_2]}e^{\sqrt{2}B_\alpha(t)-|t|^\alpha}}\in(0,\IF),
$$
where $ S_1,S_2 \in [0,\IF)$ are constants and  $B_\alpha(t),t\in\mathbb{R}$ is a standard fractional Brownian motion (fBm)  with Hurst index $\alpha/2\in(0,1].$
Further, define for $f(\cdot)$ non-negative continuous function.
\BQNY
&&\mathcal{P}_{\alpha,a}^f[-S_1,S_2]=\E{\sup_{t\in[-S_1,S_2]}e^{\sqrt{2a}B_\alpha(t)-a|t|^\alpha-f(t)}},
\EQNY
and
\BQNY
\mathcal{P}_{\alpha,a}^f[0,\IF)=\lim_{S_2\rw\IF} \mathcal{P}_{\alpha,a}^f[0,S_2],\ \
\mathcal{P}_{\alpha,a}^f(-\IF,\IF)=\lim_{S_1,S_2\rw\IF} \mathcal{P}_{\alpha,a}^f[-S_1,S_2].
\EQNY
The exact values of $\mathcal{P}_{\alpha,a}^{bt^\alpha}[0,\IF)$ are known for $\alpha=1$ and $\alpha=2$, namely,
\BQNY
\mathcal{P}_{1,a}^{bt}[0,\IF)=1+\frac{a}{b} \quad \text{and} \quad \mathcal{P}_{2,a}^{bt^2}[0,\IF)=\frac{1}{2}\LT(1+\sqrt{1+\frac{a}{b}}\RT).
\EQNY
See \cite{PicandsA,Pit72, debicki2002ruin,DI2005,DE2014,DiekerY,DEJ14,Pit20, Tabis, DM, SBK, GeneralPit16} for various properties of $\mathcal{H}_{\alpha}$ and $\mathcal{P}_{\alpha,a}^f$.
\COM{We shall write  $\mathcal{P}_{\alpha}^h$ instead of $\mathcal{P}_{\alpha}^h[0,\IF]$ and
$\mathcal{P}_{\alpha}^h$with $a>0,b\in\mathbb{R}$ and $0\leq\gamma\neq\alpha$. We shall write
 Define the generalized Piterbarg constants (if exist)
$$
\cEc{\mathcal{P}_{\alpha}^f}=\lim_{S\rightarrow\IF}\mathcal{P}_{\alpha}^f[0,S],\ \ \
\widetilde{\mathcal{P}}_{\alpha}^f =\lim_{S\rightarrow\IF}\mathcal{P}_{\alpha}^f [-S,S].
$$
}
\COM{Further, for some  closed interval $F \subset \R$ and $F\ni 0$, define regularly condition function collection $C^*_0(F)$: \\
If $f(\cdot) \in C^*_0(F)$, then $f(\cdot)$ is non-negative continuous over $\widetilde{F}$ where $\widetilde{F}\supseteq F$ is some interval
  on $\R$, $f(0)=0$, and for some $\epsilon_2>\epsilon_1>0$, $$\lim_{|t|\rw\IF,t\in F}|f(t)|/|t|^{\epsilon_1}=\IF, \lim_{|t|\rw\IF,t\in F}|f(t)|/|t|^{\epsilon_2}=0.$$}

Through this paper $\sim$ means asymptotic equivalence when the argument tends to $0$ or $\IF$.
 We notice that $\Psi(\cdot)$ denotes the tail distribution function of an $\mathcal{N}(0,1)$ random variable and  $\Psi(u)\sim\frac{1}{\sqrt{2\pi}u}e^{-\frac{u^2}{2}}, u\rw\IF$.\\
For the $L^p$ norm process $Z(t)$ in \eqref{ZZ} and a continuous function $g(t), t\in\R$, we shall investigate the asymptotics of
\BQN
\pk{\sup_{t\in[0,T]} \LT(Z^c(t)+g(t)\RT)>u},\quad u\rw\IF,
\EQN
with $c>0$ a constant.
As in \cite{FatalovA1993,Pitchi1994}, for $p\in[1,\IF]$, using the duality property of $L^p$ norm we find
\BQNY
\pk{\sup_{t\in[0,T]} Z^c(t)>u}
=\pk{\sup_{t\in[0,T]} Z(t)>u^{1/c}}
=\pk{\sup_{(t,\vk{v})\in[0,T]\times\mathcal{S}_q} Y(t,\vk{v})>u^{1/c}},
\EQNY
where $Y(t,\vk{v})=\sum_{i=1}^n d_iv_iX_{i}(t)$ is a centered Gaussian field defined on cylinder $[0,T]\times\mathcal{S}_q$ with
\BQN\label{sqq}
\mathcal{S}_q=\{\vk{v}\in\R^n:||\vk{v}||_q=1\},
\EQN
where $\frac{1}{p}+\frac{1}{q}=1$ if $q\in(1,\IF)$,  $q=\IF$ if $p=1$ and $q=1$ if $p=\IF$.
\BL\label{dd}
On $\mathcal{S}_q$ , $\sum_{i=1}^nd_i^2v_i^2$  attains its maximum $d^2$ at:\\
(i) for $p\in(2,\IF]$ at $2m$ points $\vk{v}_+^i, \vk{v}_-^i, i=1,\ldots,m,$ where $\vk{v}_+^i=(0,\ldots,0,1,0,\ldots,0)$ ($1$ stands at the i-th position), $\vk{v}_-^i=(0,\ldots,0,-1,0,\ldots,0)$ ($-1$ stands at the i-th position), $d=1$;\\
(ii) for $p=2$ at points on $\{\vk{v},\vk{v}\in\mathcal{S}_q, v_i=0,m+1\leq i\leq n\}$, $d=1$;\\
(iii) for $p\in[1,2)$ at $2^n$ points $\vk{z}$, where
$$\vk{z}=(z_1,\ldots,z_n),\ z_i=\pm(d_i/d)^{2/(q-2)}, \quad d=\LT[\sum_{i=1}^nd_i^{2p/(2-p)}\RT]^{(2-p)/2p},$$
( we take all possible $2^n$ combinations of signs "+" and "-" ),
where $z_i=\pm(d_i/d)^{0}=\pm1$.
\EL
The proof can be easily carried out by method of Lagrangian multipliers or referring to  \cite{FatalovA1993} [Lemma 3.1].\\
Next by \cite{randomChaos2015}, we have the following lemma.
\BEL \label{PPP}
 For the $L^p$ norm process $Z(t)$ in \eqref{ZZ}, if $\sigma^2(t_0)=\Var(X_i(t_0))=1,\ i=1\ldot n$ for some $t_0\in[0,\IF)$,  then we have that as $u\rw\IF$
\BQNY
\pk{Z^c(t_0)>u}
\sim \Psi\LT(\frac{u^{1/c}}{d}\RT)
\LT\{
\begin{array}{ll}
2^n(2-p)^{(1-n)/2},& \ \text{if}\ \ p\in[1,2),\\
\frac{\sqrt{2\pi}2^{\frac{(2-m)}{2}}u^{\frac{m-1}{c}}}
{\Gamma(m/2)}\prod_{i=m+1}^n(1-d_i^2)^{-\frac{1}{2}},& \ \text{if}\ \ p=2,\\
2m,& \ \text{if}\ \ p\in(2,\IF],
\end{array}
\RT.
\EQNY
\COM{and the density function satisfying
\BQNY
p_{||\vk{X}*\vk{d}||_p^c}(u)
\sim \frac{u^{2/c-1}}{c d^2}\Psi\LT(\frac{u^{1/c}}{d}\RT)
\LT\{
\begin{array}{ll}
2^n(2-p)^{(1-n)/2},& \ \text{if}\ \ p\in(1,2),\\
\frac{\sqrt{2\pi}2^{\frac{(2-m)}{2}}u^{\frac{m-3}{c}}}
{\Gamma(m/2)}\prod_{i=m+1}^n(1-d_i^2)^{-\frac{1}{2}},& \ \text{if}\ \ p=2,\\
2m,& \ \text{if}\ \ p\in(2,\IF),
\end{array}
\RT.
\EQNY}
with the convention $\prod_{i=n+1}^n(1-d_i^2)^{-\frac{1}{2}}=1$ and $d$ the same as in \nelem{dd}.
\EEL

\section{Extremes of $L^p$ norm processes with trend}
In this section, recall that $Z(t)$ in \eqref{ZZ} is the \emph{$L^p$ norm process} and $X_{i}(t)$'s are independent copies of $X(t)$ with continuous trajectories,  variance functions $\sigma^2(\cdot)$ and correlation functions $r(\cdot,\cdot)$.
\subsection{Extremes of non-stationary $L^p$ norm processes with trend}
As in  \cite{Threshold2016}, if $X(t)$ is non-stationary, we introduce the following assumptions:
\begin{itemize}
\item [(i)] $\sigma(\cdot)$ attains its maximum  on $[0, T]$ at the unique point $t_0\in[0,T]$ and
\BQNY
\sigma(t)=1-b|t-t_0|^{\beta}+o(|t-t_0|^{\beta}),\ \
\ \ t\rw t_0
\EQNY
for some positive constants $b,\beta$.
\item [(ii)]$
r(s,t)=1-a|t-s|^{\alpha}+o(|t-s|^{\alpha}), \ s, t \rw t_0
$
for some constants $a>0$ and $\alpha\in(0,2].$
\end{itemize}
Further, we introduce a bounded measurable trend function $g(t)$ which satisfies
\begin{itemize}
\item[(iii)]
$
g(t)\sim -w|t-t_0|^{\gamma},\ \ \ \  t \rw t_0
$
for some constants $\gamma>0, w\geq0$.
\end{itemize}
\BT\label{Thm2}
If assumptions (i)-(iii) are satisfied,
then for $\vk{d}$ in \eqref{vd} and $d$ in \nelem{dd}, we have as $u\rw\IF$
\BQNY\label{ei2}
\pk{\sup_{t\in[0,T]} \LT(Z^c(t)+g(t)\RT)>u}\sim\pk{Z^c(t_0)>u}
\LT\{
\begin{array}{ll}
u^{\frac{2}{\alpha^*}-\frac{2}{\beta^*}}a^{1/\alpha} d^{-2/\alpha}\mathcal{H}_{\alpha}\int_{Q}^\IF e^{-f(t)}dt,& \ \text{if}\ \ \alpha^*<\beta^*,\\
\mathcal{P}_{\alpha,a d^{-2}}^{f(t)}[Q,\IF),& \ \text{if}\ \ \alpha^*=\beta^*,\\
1,& \ \text{if}\ \ \alpha^*>\beta^*,\\
\end{array}
\RT.
\EQNY
where $\alpha^*=\alpha c$, $\beta^*=\min(\beta c,\frac{2\gamma c}{2-c})\II_{\{c<2\}}+\beta c\II_{\{c\geq 2\}}$, $f(t)=\frac{b|t|^\beta}{d^2}\II_{\{\beta^*=\beta c\}}+\frac{w|t|^{\gamma}}{c d^2}\II_{\{\beta^*=\frac{2\gamma c}{2-c}\}}$, and $ Q= -\IF$ if $t_0\in(0,T)$, $ Q=0$ if $t_0\in\{0,T\}$.
\ET
\begin{remarks}\label{rem11}
i) In \netheo{Thm2}, if we assume that $w=0$, we get the extremes of centered non-stationary $L^p$ norm processes i.e.,
\BQNY
\pk{\sup_{t\in[0,T]} Z^c(t)>u},\ u\rw\IF.
\EQNY
ii) Following the similar arguments as in the proof of \netheo{Thm2}, the result in \netheo{Thm2} still holds for $w<0, \gamma\geq \frac{(2-c)\beta}{2}$ if $c<2$ and $w<0, \gamma>0$ if $c\geq 2$.
\end{remarks}
%concerning the asymptotics of $\phi^{p}_{c,T}(u)$ as $u\to\IF.$
%where $\vk{X}(t)$ is a centered vector Gaussian process, $g(t)$ is a continuous trend function and $\vk{d}$ is the same as in \eqref{vd}.
\COM{In order to analyze the asmptotics of \eqref{eq:Xgt1}, define below for large enough $u$
\BQN\label{eq:Mtu}
\sigma_{u}(t):=\frac{\sigma(t)}{(1-{g(t)}/{u})^{1/c}},\ \ \ t\in[0,S].
\EQN}
%\subsection{Locally-stationary case}
\subsection{Extremes of locally stationary $L^p$ norm processes with trend}
If $X(t)$ is locally stationary, as in \cite{Threshold2016}, we shall suppose that:
\begin{itemize}
\item[(iv)]
$
r(s,t)=1-a(t)|t-s|^{\alpha} +o(|t-s|^{\alpha}),\ \ \abs{t-s}\rw 0,\ \alpha\in(0,2],
$
where $a(t)$ are positive continuous function on $[0,T]$.
\item[(v)] $
r(s,t)<1,\ \  \forall\ s,t\in[0,T] \quad \hbox{and}\quad s\neq t.
$
\end{itemize}
Before giving the scenarios with trend, we consider the extremes of the centered locally stationary $L^p$ norm processes.
\BT\label{Thm0} Assume that $\sigma(t)\equiv1$, i.e., unit variance and covariance function $r(\cdot,\cdot)$ satisfies assumptions (iv) and (v).
Then we have for $c>0$
\BQNY
\pk{\sup_{t\in[0,T]} Z^c(t)>u}
\sim \int_{0}^T(a(t))^{\frac{1}{\alpha}}dt d^{-\frac{2}{\alpha}} \mathcal{H}_\alpha u^{\frac{2}{\alpha c}}\pk{Z^c(0)>u},\ u\rw \IF,
\EQNY
where $d$ is the same as in \nelem{dd}.
\ET
\BT\label{Thm1}
Assume that $\sigma(t)\equiv1$, i.e., unit variance and correlation function $r(\cdot,\cdot)$ satisfies assumptions (iv) and (v). Assume that $g(t)\ t\in[0,T]$ is a continuous function which attains its maximum at a unique point  $t_0\in[0,S]$ satisfying assumption (iii) for some constants $w, \gamma>0$. Further, set $\alpha^*=\alpha c,\ \beta^*=\frac{2\gamma c}{2-c}\II_{\{c<2\}}$ and  $f(t)=\frac{w \abs{t}^\gamma}{c d^{2}}$ and $d$ is the same as in \nelem{dd}.\\
If $c\in(0,2)$, then we have as $u\rw\IF$
\BQNY
\pk{\sup_{t\in[0,T]} \LT(Z^c(t)+g(t)\RT)>u}\sim u^{(\frac{2}{\alpha^*}-\frac{2}{\beta^*})_{+}}\pk{Z^c(0)>u}
\LT\{
\begin{array}{ll}
a^{\frac{1}{\alpha}}d^{-\frac{2}{\alpha}}\mathcal{H}_{\alpha}\int_{Q}^\IF e^{-f(t)}dt,& \text{if}\ \alpha^*<\beta^*,\\
\mathcal{P}_{\alpha,a d^{-2}}^{f(t)}[Q,\IF), & \text{if}\ \alpha^*=\beta^*,\\
1,& \text{if}\ \alpha^*>\beta^*,
\end{array}
\RT.
\EQNY
where $a=a(t_0)$ and $ Q= -\IF$ if $t_0\in(0,T)$, $ Q=0$ if $t_0\in\{0,T\}$.\\
If $c=2$, then we have
\BQNY
\pk{\sup_{t\in[0,T]} \LT(Z^c(t)+g(t)\RT)>u}\sim
\int_{0}^T(a(t))^{\frac{1}{\alpha}}e^{\frac{g(t)}{2d^{2}}}dt d^{-\frac{2}{\alpha}} \mathcal{H}_\alpha u^{\frac{2}{\alpha^*}}\pk{Z^c(0)>u},\ u\rw \IF.
\EQNY
If $c>2$, then we have
\BQNY
\pk{\sup_{t\in[0,T]} \LT(Z^c(t)+g(t)\RT)>u}\sim
\int_{0}^T(a(t))^{\frac{1}{\alpha}}dt d^{-\frac{2}{\alpha}} \mathcal{H}_\alpha u^{\frac{2}{\alpha^*}}\pk{Z^c(0)>u},\ u\rw \IF.
\EQNY
\ET
\begin{remark}
By the proof, we notice that for the case $c=2$ in \netheo{Thm1}, the result always holds for any continuous function $g(t), t\in[0,1]$. When $c>0$, the result holds for any bounded function $g(t), t\in[0,1]$.
\end{remark}

\begin{example}\label{EX1}
For $Z(t)$ in \eqref{ZZ} with $X_i(t)=B^i_\alpha(t),i=1\ldot n$ the independent fractional Brownian motions, we have as $u\rw\IF$
\BQNY
\pk{\sup_{t\in[0,1]} Z(t)-\sqrt{1-t}>u}=\pk{ Z(1)>u}
\LT\{
\begin{array}{ll}
u^{\frac{2}{\alpha}-2}\LT(\frac{1}{2d^2}\RT)^{1/\alpha} \mathcal{H}_{\alpha}\int_{0}^\IF e^{-f(t)}dt,& \ \text{if}\ \ \alpha<1,\\
\mathcal{P}_{\alpha, d^{-2}/2}^{f(t)}[0,\IF),& \ \text{if}\ \ \alpha=1,\\
1,& \ \text{if}\ \ \alpha>1,\\
\end{array}
\RT.
\EQNY
where $f(t)=\frac{\alpha}{2d^2}t+\frac{1}{d^2}t^{\frac{1}{2}}$ and $d$ is the same as in \nelem{dd}
\end{example}
Following example is a special case of \netheo{Thm1}, which is corresponded with \cite{LJ2017} [Theorem 2.1].
\begin{example}\label{EX2}
In \netheo{Thm1}, assume that $p=2$, $c=2$ and $g(t), t\in[0,T]$ is a continuous function, then we have
\BQNY
\pk{\sup_{t\in[0,T]} \LT(Z_2^2(t)+g(t)\RT)>u}\sim
\int_{0}^T(a(t))^{\frac{1}{\alpha}}e^{\frac{g(t)}{2}}dt \mathcal{H}_\alpha \frac{2^{1-m/2}\prod_{i=m+1}^n(1-d_i^2)^{-\frac{1}{2}}}
{\Gamma(m/2)}u^{\frac{m-1}{2}+\frac{1}{\alpha}}e^{-\frac{u}{2}},\ u\rw \IF.
\EQNY
\end{example}

\section{Applications}

\subsection{Ruin probability of a risk model}\label{ruin}
In theoretical insurance modelling a surplus process $U(t)$ can be defined by
\BQNY
U(t)=u+wt-X(t), \ \ \ t\ge0,
\EQNY
see \cite{MR1458613}, where $u\geq 0$ is the initial reserve, $w>0$ is the rate of premium and the stochastic process $ X(t),t\geq0 $ denotes the aggregate claims process. See \cite{rolski2009stochastic, DHJ13a,HX2007,ParisianBrownianfinite2017,Threshold2016,ParisianInfinite2018} for more studies on related risk models.
Here we investigate   
$$X(t)=\sum_{i=1}^n \abs{d_iB_\alpha^{i}(t)}^2, \quad t\ge0,$$
where $\vk{d}=(d_1\ldot d_n)$ is the same as in \eqref{vd} and $B_\alpha^{i}(t)$ are independent fractional Brownian motions. $X(t)$ can be considered as the sum of $n$ independent claims or payments until time $t$.
The corresponding ruin probability over a finite-time horizon $[0,1]$ is defined as
$$\pk{\inf_{t\in[0,1]}U(t)<0}. $$
We present next approximation of this ruin probability. 
\begin{prop}\label{riskmodel}
We have as $u\rw\IF$
\BQNY
\pk{\inf_{t\in[0,1]}U(t)<0}\sim u^{m/2-1}e^{-\frac{u+w}{2}}\frac{2^{1-m/2}}{\Gamma(m/2)}\prod_{i=m+1}^n\LT(1-d_i^2\RT)^{-1/2}\LT\{
\begin{array}{ll}
u^{1/\alpha-1}2^{1-1/\alpha} \mathcal{H}_{\alpha},& \ \text{if}\ \ \alpha<1,\\
2,& \ \text{if}\ \ \alpha=1,\\
1,& \ \text{if}\ \ \alpha>1.\\
\end{array}
\RT.
\EQNY
\end{prop}

Besides in risk modelling, the $L^p$ norm processes, especially the chi-square processes, are also widely utilized in hypothesis testing, see \cite {HTDa1987,SCSP2014} and the reference there. Next we give an example.

\subsection{The Ornstein-Uhlenbeck chi-square process in Quantitative Trait Locus detection}
\COM{In \cite{}, the focus is on hypothesis testing when a nuisance parameter $t^*$ is present only under an alternative. So, $t^*$ is meaningless under the null hypothesis. If $t^*$ were known, the natural way to perform the test is to consider $t=t^*$. However, as it is only known that $t$ belongs to the interval $[S_1,S_2]$, \cite{} suggests the use of the test statistic:
\BQNY
\{S(t): S_1\leq t<\leq S_2 \}
\EQNY
where $S_1<S_2$ are two constants and $S(t)$ can be consider as 
\BQNY
S(t)=\sum_{i=1}^{n} B_i(t)^2
\EQNY
with $B_i(t)$ are independent Brownian motions.}
A Quantitative Trait Locus (QTL) denotes a gene with quantitative effect on a trait. The method used by most of geneticists in order to detect a QTL on a chromosome, is the Interval Mapping proposed by \cite{LanderBO1989}. Using the Haldane distance and modelling in \cite{Haldane1919}, each chromosome is represented by a segment $[0,T]$. The distance on $[0,T]$ is called the genetic distance. At each location $t\in[0,T]$, using the "genome information" brought by genetic markers, a likelihood ratio test (LRT) is performed, testing the presence of a QTL at this position. \cite{ADR2014} prove that when the number of genetic markers and the number of progenies tends to infinity, the limiting process of the LRT process is an Ornstein-Uhlenbeck chi-square process under the null hypothesis of the absence of QTL on the interval $[0,T]$. In order to take decision about the presence of a QTL on $[0,T]$, we need to calculate the critical value for the supremum of an Ornstein-Uhlenbeck chi-square process, i.e.,
\BQNY
\sup_{t\in[0,T]}S(t),
\EQNY
where the Ornstein-Uhlenbeck chi-square process $S(t)$ is 
\BQNY
S(t)=\sum_{i=1}^{n} V_i(t)^2
\EQNY
and $V_i(t),\ 1\leq i\leq n$ are independent identically stationary Gaussian processes with covariance function given by 
$$Cov(V_i(s),V_i(t))=e^{-2\abs{t-s}}.$$
\begin{prop}\label{OUCSP}
We have as $u\rw\IF$
\BQNY
\pk{\sup_{t\in[0,T]}S(t)>u}\sim\frac{2^{2-n/2}}{\Gamma(n/2)}T u^{n/2}e^{-u/2} .
\EQNY
\end{prop}

\section{Proofs}

During the following proofs, $\mathbb{Q}_i, i\in\N$ are some positive constants which can be different from line by line and for  interval $\Delta_1,\Delta_2  \subseteq[0,\IF)$ we denote
\BQNY
&&\mathcal{L}_u(\Delta_1):=\pk{\sup_{t\in\Delta_1} \LT(Z^c(t)+g(t)\RT)>u}, \ \
\mathcal{L}_u(\Delta_1,\Delta_2):=\pk{\sup_{t\in\Delta_1} \LT(Z^c(t)+g(t)\RT)>u,\sup_{t\in\Delta_2} \LT(Z^c(t)+g(t)\RT)>u},
\EQNY
and
\BQNY
\mathcal{K}_u(\Delta_1):=\pk{\sup_{t\in\Delta_1}Z^c(t)>u}, \ \
\mathcal{K}_u(\Delta_1,\Delta_2):=\pk{\sup_{t\in\Delta_1}Z^c(t)>u,\sup_{t\in\Delta_2} Z^c(t)>u}.
\EQNY
\prooftheo{Thm2}
We first present the proof for the case  $t_0=0$.\\
Set $\beta^*=\min(\beta c,\frac{2\gamma c}{2-c})\II_{\{c<2\}}+\beta c\II_{\{c\geq 2\}},\ \alpha^*=\alpha c$, $\delta(u)=\frac{(\ln u)^\rho }{u^{2/\beta^*}}$ with $\rho>\max\LT(\frac{1}{\beta},\frac{1}{\gamma}\RT)$ and for $u$ large enough
$$Y(t,\vk{v})=\sum_{i=1}^n d_iv_iX_{i}(t),\quad (t,\vk{v})\in\R\times\mathcal{S}_q$$
with $\mathcal{S}_q$ the same as in \eqref{sqq} which is  a centered Gaussian field.\\
We have for some small $\theta>0$ and $u$ large enough
\BQN\label{bb1}
\mathcal{L}_u([0,\delta(u)])\leq\mathcal{L}_u([0,T])\leq \mathcal{L}_u([0,\delta(u)])+\mathcal{L}_u([\delta(u),\theta])+
\mathcal{L}_u([\theta,T]).
\EQN
We  first give the upper bounds of $\mathcal{L}_u([\delta(u),\theta])$ and $\mathcal{L}_u([\theta,T])$.\\
Set $\sigma_\theta:=\sup_{t\in[\theta,T]}\sigma(t)<1$ and $g_m=\sup_{t\in[0,T]}g(t)<\IF$.
Then by Borell inequality as in \cite{AdlerTaylor} and \nelem{PPP} for large $u$
\BQN\label{PI2}
\mathcal{L}_u([\theta,T])&\leq&\pk{\sup_{t\in[\theta,T]}Z(t)>(u-g_m)^{1/c}}\nonumber\\
&\leq& \pk{\sup_{(t,\vk{v})\in [\theta,T]\times\mathcal{S}_q} Y(t,\vk{v})>(u-g_m)^{1/c}}\nonumber\\
&\leq& \exp\LT(-\frac{\LT((u-g_m)^{1/c}-\mathbb{Q}_1\RT)^2}{2V^*_Y}\RT)\nonumber\\
&=&o\LT(\pk{Z^c(0)>u}\RT),\ u\rw\IF,
\EQN
where $\mathbb{Q}_1:=\E{\sup_{(t,\vk{v})\in [\theta,T]\times\mathcal{S}_q}
Y(t,\vk{v})}<\IF$ and
$$V^*_Y:=\sup_{(t,\vk{v})\in [\theta,T]\times\mathcal{S}_q}\Var \LT(Y(t,\vk{v})\RT)\leq \LT(\sup_{t\in[\theta,T]}\sigma^2(t)\RT)d^2=\sigma^2_\theta d^2<d^2.$$
By assumptions (i) and (iii), we know that for some $\vn_1\in(0,1)$
\BQN
&\frac{u-g(t)}{\sigma^c(t)}
\geq (u+w(1-\vn_1)|t|^\gamma)(1+(1-\vn_1)bc|t|^\beta)
\geq u\LT(1+\frac{w(1-\vn_1)}{u}|t|^\gamma+(1-\vn_1) bc|t|^\beta\RT),\label{hl2}\\
&\frac{u-g(t)}{\sigma^c(t)}
\leq (u+w(1+\vn_1)|t|^\gamma)(1+(1+\vn_1)bc|t|^\beta)
\leq u\LT(1+\frac{w(1+\vn_1)}{u}|t|^\gamma+(1+\vn_1) bc|t|^\beta\RT)\label{hu2}
\EQN
hold for $t\in[0,\theta]$ when $\theta$ small enough, then
\begin{align}\label{boundsig1}
\inf_{t\in[\delta(u),\theta] }\frac{(u-g(t))^{2/c}}{\sigma^2(t)}
&\geq \inf_{t\in[\delta(u),\theta] } u^{2/c}\LT(1+\frac{w(1-\vn_1)}{u}|t|^\gamma+(1-\vn_1) bc|t|^\beta\RT)^{2/c}\nonumber\\
&\geq u^{2/c}+\mathbb{Q}_2(\ln u)^{(\rho \beta)\vee(\rho \gamma)}.
\end{align}

Denote $\overline{\X}(t)=(\overline{X}_1(t)\ldot \overline{X}_n(t))$ with $\overline{X}_i(t)=\frac{X_i(t)}{\sigma(t)},\  t\in[0,\theta]$. By assumption (ii), we have that
\begin{align*}
\E{\LT(\LT(\sum_{i=1}^n d_iv_i\overline{X}_{i}(t)\RT)-\LT(\sum_{i=1}^n d_iv'_i\overline{X}_{i}(s)\RT)\RT)^2}
&\leq 2\E{\LT(\LT(\sum_{i=1}^n d_iv_i\overline{X}_{i}(t)\RT)-\LT(\sum_{i=1}^n d_iv_i\overline{X}_{i}(s)\RT)\RT)^2}\\
&\quad+2\E{\LT(\LT(\sum_{i=1}^n d_iv_i\overline{X}_{i}(s)\RT)-\LT(\sum_{i=1}^n d_iv'_i\overline{X}_{i}(s)\RT)\RT)^2}\\
&\leq 4\E{\sum_{i=1}^n\LT(\overline{X}_{i}(t)-\overline{X}_{i}(s)\RT)^2}
+4\E{\sum_{i=1}^n (v_i-v'_i)^2\LT(\overline{X}_{i}(s)\RT)^2}\\
&\leq \mathbb{Q}_3|s-t|^\alpha+ \mathbb{Q}_4\sum_{i=1}^n|v_i-v'_i|^2\\
&\leq \mathbb{Q}_5\LT(|s-t|^\alpha+\sum_{i=1}^n|v_i-v'_i|^\alpha\RT)
\end{align*}
holds for $s,t\in[0,\theta]$ and $\vk{v},\vk{v}'\in\mathcal{S}_q$.
Thus it follows from \cite{Pit96} [Theorem 8.1], \eqref{boundsig1} and \nelem{PPP} that
\begin{align}\label{PI3}
\mathcal{L}_u([\delta(u),\theta])
&\leq \pk{\sup_{(t,\vk{v})\in[\delta(u),\theta]\times\mathcal{S}_q} \sum_{i=1}^n d_iv_i\overline{X}_{i}(t)>\inf_{s\in[\delta(u),\theta]}\frac{(u-g(s))^{1/c}}{\sigma(s)}}
\nonumber\\
&\leq\mathbb{Q}_6u^{\frac{2(n+1)}{\alpha}}\Psi\LT(\inf_{s\in[\delta(u),\theta]}
\frac{(u-g(s))^{1/c}}{d\sigma(s)}\RT)\nonumber\\
&\leq \frac{\mathbb{Q}_6}{d\sqrt{2\pi}}u^{\frac{2(n+1)}{\alpha}-\frac{2}{c}}
\exp\LT(-\frac{1}{2d}\LT( u^{2/c}+\mathbb{Q}_2(\ln u)^{(\rho \beta)\vee(\rho \gamma)}\RT)\RT)\nonumber\\
&=o\LT(\pk{Z^c(0)>u}\RT),\ u\rw\IF.
\end{align}
Thus by \eqref{PI2}, \eqref{PI3} and the fact that $\mathcal{L}_u([0,\delta(u)])\geq \pk{Z^c(0)>u}$ for $u$ positive, we have
\BQN\label{Thmeq02}
\mathcal{L}_u([\delta(u),\theta])=o\LT(\mathcal{L}_u([0,\delta(u)])\RT),\ \ \mathcal{L}_u([\theta,T])=o\LT(\mathcal{L}_u([0,\delta(u)])\RT),\ u\rw\IF ,
\EQN
which combined with \eqref{bb1} imply
\BQN\label{Main1}
\mathcal{L}_u([0,T])\sim \mathcal{L}_u([0,\delta(u)]),\ u\rw\IF.
\EQN

Now we focus on the asymptotic of $\mathcal{L}_u([0,\delta(u)])$, as $u\rw\IF$.\\
\COM{When $p\in(1,2)\cup(2,\IF)$, by \nelem{dd}, we know $\sigma_u^2(t,\vk{v}):=\Var\LT(Y^u(t,\vk{v})\RT)=\frac{\sigma^2(t)}{(1-g(t)/u)^{2/c}}\LT(\sum_{i=1}^n d_i^2v_i^2\RT)$ attains the maximum over $E_u$ at several discrete points, so we can choose $\vn$ small enough such that
$E_u^{\vn}(i)=[0,\theta]\times\mathcal{S}^\vn_q(i) $ with $\mathcal{S}^\vn_q(i)$ the union of non-overlapping compact neighborhoods of $ \vk{v}_+^i, \vk{v}_-^i$ or $\vk{z}$ in \nelem{dd}.
Then as mentioned in \cite{Pit96} or \cite{Fatalov1992}[Lemma 2.1]
\BQN\label{sum1}
\pk{\sup_{(t,\vk{v})\in [0,\theta]\times\mathcal{S}^{\vn}_q}Y^u(t ,\vk{v})>u}\sim
\sum_{i=1}^M\pk{\sup_{(t,\vk{v})\in [0,\theta]\times\mathcal{S}^{\vn}_q(i)}Y^u(t ,\vk{v})>u},
\EQN
where $M$ is the number of the maximum points of $\sigma_u^2(t,\vk{v})$.\\
{\bf Step 1:} $p\in(1,2)$ and $M=2^k$.\\
It is enough to find the asymptotics of single term in \eqref{sum1}, for instance, for a point $(0,\vk{z}),\ z_i=(d_i/d)^{2/q-2}$. In a neighborhood $\mathcal{S}^\vn_q(1)$ of $\vk{z}$, we have
$$v_n=\LT(1-\sum_{i=1}^{n-1}v_i^q\RT)^{1/q},$$
hence the fields $\frac{Y(u^{-2/\alpha}t+t_0 ,\vk{v})}{1+u^{-2}f(t)}$ can be represented as
\BQNY
Y_1(u^{-2/\alpha}t+t_0, \widetilde{\vk{v}})
=\sum_{i=1}^{n-1}v_id_i \frac{X_i(u^{-2/\alpha}t+t_0)}{1+u^{-2}f(t)}
+\LT(1-\sum_{i=1}^{n-1}v_i^q\RT)^{1/q}d_n\frac{X_n(u^{-2/\alpha}t+t_0)}{1+u^{-2}f(t)},
 \widetilde{\vk{v}}=(v_1,\cdots,v_{n-1}),
\EQNY
which is defined in $[-S_1,S_1]\times\widetilde{\mathcal{S}}^\delta_q(1)$ where
$$\widetilde{\mathcal{S}}^\delta_q(1)
=\LT\{\widetilde{\vk{v}}:\LT(v_1,\cdots,v_{n-1},\LT(1-\sum_{i=1}^{n-1}v_i^q\RT)^{1/q}\RT)\in\mathcal{S}^\delta_q(1)\RT\},$$
is a small neighborhood of $\widetilde{\vk{z}}=\LT(z_1,\cdots,z_{n-1}\RT)$.
On $[-S_1,S_1]\times\widetilde{\mathcal{S}}^\delta_q(1)$, the variance
\BQNY
\sigma_1^2(t,\widetilde{\vk{v}}):=\frac{1}{(1+u^{-2}f(t))^2}\sigma_1^2(\widetilde{\vk{v}})
:=\frac{1}{(1+u^{-2}f(t))^2}\LT[\sum_{i=1}^{n-1}d_i^2v_i^2+d_n^2\LT(1-\sum_{i=1}^{n-1}v_i^q\RT)^{2/q}\RT]
\EQNY
of $Y_1(u^{-2/\alpha}t+t_0, \widetilde{\vk{v}})$ attains its maximum $d^2$ at $(0,\widetilde{\vk{z}})$ where $\widetilde{\vk{z}}$ is a interior point of a set $\widetilde{\mathcal{S}}^\delta_q(1)$. We can write the following Taylor expansion for $\sigma_1(t,\widetilde{\vk{v}})$
\BQNY
\sigma_1(t,\widetilde{\vk{v}})=\frac{d}{1+u^{-2}f(t)}
-\frac{q-2}{2d}(\widetilde{\vk{v}}-\widetilde{\vk{z}})
\Lambda(\widetilde{\vk{v}}-\widetilde{\vk{z}})^T+o(|\widetilde{\vk{v}}-\widetilde{\vk{z}}|^2), \widetilde{\vk{v}}\rw\widetilde{\vk{z}},\ u\rw\IF,
\EQNY
where $\Lambda=(\lambda_{i,j})_{i,j=1,\cdots,n-1}$ is a non-negative define matrix with elements
\BQNY
\lambda_{i,j}=-(2(q-2))^{-1}\frac{\partial^2}{\partial v_i\partial v_j}\LT[\sum_{i=1}^{n-1}d_i^2v_i^2+d_n^2\LT(1-\sum_{i=1}^{n-1}v_i^q\RT)^{2/q}\RT]|_{\widetilde{\vk{v}}=\widetilde{\vk{z}}},
i,j=1,\cdots,n-1.
\EQNY
We have the following expansion for the correlation function $r_1(t,\widetilde{\vk{v}},s,\widetilde{\vk{v}}')$ of $Y_1(u^{-2/\alpha}t+t_0, \widetilde{\vk{v}})$
\BQNY
r_1(t,\widetilde{\vk{v}},s,\widetilde{\vk{v}}')=1- u^{-2}a(t-s)^\alpha
-\frac{1}{2d}(\widetilde{\vk{v}}-\widetilde{\vk{z}})
\Lambda(\widetilde{\vk{v}}-\widetilde{\vk{z}})^T+o(|\widetilde{\vk{v}}-\widetilde{\vk{z}}|^2), \widetilde{\vk{v}},\widetilde{\vk{v}}'\rw\widetilde{\vk{z}}, u\rw\IF.
\EQNY
There exists a non-singular matrix $Q$ such that $Q\Lambda Q^T$ is diagonal, and
set the diagonal is $(c_1,\cdots,c_{n-1})$. Then
\BQNY
\sigma_1(t,Q\widetilde{\vk{v}})=
d-du^{-2}f(t)-\frac{q-2}{2d}\sum_{i=1}^{n-1}c_i(v_i-z_i)^2
+o(|\widetilde{\vk{v}}-\widetilde{\vk{z}}|^2), \widetilde{\vk{v}}\rw\widetilde{\vk{z}},\ u\rw\IF,
\EQNY
and
\BQNY
r_1(t,Q\widetilde{\vk{v}},s,Q\widetilde{\vk{v}}')=1- u^{-2}a(t-s)^\alpha
-\frac{1}{2d}\sum_{i=1}^{n-1}c_i(v_i-z_i)^2+o(|\widetilde{\vk{v}}-\widetilde{\vk{z}}|^2), \widetilde{\vk{v}},\widetilde{\vk{v}}'\rw\widetilde{\vk{z}}, u\rw\IF.
\EQNY
Then set $Y_2(u^{-2/\alpha}t+t_0,\widetilde{\vk{v}})=Y_1(u^{-2/\alpha}t+t_0,Q\widetilde{\vk{v}})$,
defined on a set $[-S_1,S_2]\times(Q^{-1}\widetilde{\mathcal{S}}^\delta_q(1))$. We know that the point $Q\widetilde{\vk{z}}$ is a interior point of $Q^{-1}\widetilde{\mathcal{S}}^\delta_q(1)$.
Then the proof follows by similar arguments as in the proof of \cite{Pit96} [Theorem 8.2]. Consequently, we get
\BQNY
\pk{\sup_{(t,\vk{v})\in \mathcal{D}^1_\delta}\frac{Y(u^{-2/\alpha}t+t_0 ,\vk{v})}{1+u^{-2}f(t)}>u}
&=&\pk{\sup_{(t,\widetilde{\vk{v}})\in [-S_1,S_2]\times(Q^{-1}\widetilde{\mathcal{S}}^\delta_q(1))}Y_2(u^{-2/\alpha}t+t_0, \widetilde{\vk{v}})>u}\\
&\sim&\mathcal{P}_{\alpha,ad^{-2}}^{\frac{1}{d^{2}}f(t)}
[-S_1,S_2]\LT(\prod_{i=1}^{n-1}\mathcal{P}_{2,1}^{(q-2)t^2}
(-\IF,\IF)\RT)\Psi\LT(\frac{u}{d}\RT)\\
&=&\mathcal{P}_{\alpha,ad^{-2}}^{\frac{1}{d^{2}}f(t)}
[-S_1,S_2](2-p)^{(1-n)/2}\Psi\LT(\frac{u}{d}\RT),
\EQNY
 and
 \BQNY
 \pk{\sup_{(t,\vk{v})\in [-S_1,S_2]\times\mathcal{S}_q}\frac{Y(u^{-2/\alpha}t+t_0 ,\vk{v})}{1+u^{-2}f(t)}>u}\sim 2^n\mathcal{P}_{\alpha,ad^{-2}}^{\frac{1}{d^{2}}f(t)}
[-S_1,S_2](2-p)^{(1-n)/2}\Psi\LT(\frac{u}{d}\RT).
 \EQNY}
Denote for any $\lambda>0$ and some $\vn\in(0,1)$
\BQNY
&&I_k(u)=[ku^{-2/\alpha^*}\lambda,(k+1)u^{-2/\alpha^*}\lambda],\ \  k\in\N,\ \  N(u)=\LT\lfloor(\ln u)^{\frac{2q}{\beta^*}}u^{\frac{2}{\alpha^*}-\frac{2}{\beta^*}}\lambda^{-1}\RT\rfloor,\\
&&\mathcal{G}_{u,+\vn}(k)=u\LT(1+\frac{w(1+\vn)}{u}\abs{(k+1)u^{-2/\alpha^*}\lambda}^\gamma+(1+\vn) bc\abs{(k+1)u^{-2/\alpha^*}\lambda}^\beta\RT), \\
&&\mathcal{G}_{u,-\vn}(k)=u\LT(1+\frac{w(1-\vn)}{u}\abs{ku^{-2/\alpha^*}\lambda}^\gamma+(1-\vn) bc\abs{ku^{-2/\alpha^*}\lambda}^\beta\RT).
\EQNY
{\bf Case 1}: $\beta^*>\alpha^*$. For $u$ large enough, we have
\BQN
\sum_{k=0}^{N(u)-1}\mathcal{L}_u(I_{k}(u))
-\sum_{i=1}^2\Lambda_i(u)\leq\mathcal{L}_u([0,\delta(u)])
\leq\sum_{k=0}^{N(u)}\mathcal{L}_u(I_{k}(u)),\label{Thmeq12}
\EQN
where
$$\Lambda_1(u)=\sum_{k=0}^{N(u)}\mathcal{L}_u(I_{k}(u),I_{k+1}(u)),\quad
\Lambda_2(u)=\sum_{0\leq k,l\leq N(u), l\geq k+2}\mathcal{L}_u(I_{k}(u),I_{l}(u)).$$
In the view of \nelem{im1} and \eqref{hl2}, we have that for some $\epsilon\in [0,1)$,
\BQN\label{Th21}
\sum_{k=0}^{N(u)}\mathcal{L}_u(I_{k}(u))
&\leq&\sum_{k=0}^{N(u)}\pk{\sup_{t\in I_k(u)}\LT\|\overline{\vk{X}}(t)*\vk{d}\RT\|_p^c>\mathcal{G}_{u,-\vn}(k)}\nonumber\\
&\sim& \mathcal{H}_\alpha[0,a^{1/\alpha}d^{-2/\alpha}\lambda]
\sum_{k=0}^{N(u)}\pk{Z^c(0)>\mathcal{G}_{u,-\vn}(k)}\nonumber\\
&\sim& \mathcal{H}_\alpha[0,a^{1/\alpha}d^{-2/\alpha}\lambda]\pk{Z^c(0)>u}
\nonumber\\
&&\times
\sum_{k=0}^{N(u)}\exp\LT(-(1-\vn-\epsilon)\frac{w}{cd^2}u^{\frac{2-c}{c}}|k S u^{-\frac{2}{\alpha^*}}|^\gamma -(1-\vn-\epsilon)\frac{b}{d^2}u^{2/c}|k S u^{-\frac{2}{\alpha^*}}|^\beta\RT)\nonumber\\
&\sim& \mathcal{H}_\alpha[0,a^{1/\alpha}d^{-2/\alpha}\lambda]\pk{Z^c(0)>u}
\sum_{k=0}^{N(u)}\exp\LT(-(1-\vn-\epsilon)f(u^{\frac{2}{\beta^*}}k S u^{-\frac{2}{\alpha^*}})\RT)\nonumber\\
&\sim&
\pk{Z^c(0)>u}\frac{\mathcal{H}_\alpha[0,a^{1/\alpha}d^{-2/\alpha}\lambda]}{\lambda}u^{\frac{2}{\alpha^*}-\frac{2}{\beta^*}}
\int_0^{\IF}\exp\LT(-(1-\vn-\epsilon)f(t)\RT)dt\nonumber\\
&\sim&
\pk{Z^c(0)>u}a^{1/\alpha}d^{-2/\alpha}\mathcal{H}_\alpha u^{\frac{2}{\alpha^*}-\frac{2}{\beta^*}}
\int_0^{\IF}e^{-f(t)}dt,
\EQN
as $ u\rw\IF, \ \lambda\rw\IF,\ \vn\rw 0, \ \epsilon\rw 0$ where
$f(t)=\frac{b|t|^\beta}{d^2}\II_{\{\beta^*=\beta c\}}+\frac{w}{c d^2}|t|^{\gamma}\II_{\{\beta^*=\frac{2\gamma c}{2-c}\}}$. Similarly, we derive that
\BQN\label{Th212}
\sum_{k=0}^{N(u)-1}\mathcal{L}_u(I_{k}(u)) \geq\pk{Z^c(0)>u}a^{1/\alpha}d^{-2/\alpha}\mathcal{H}_\alpha u^{\frac{2}{\alpha^*}-\frac{2}{\beta^*}}
\int_0^{\IF}e^{-f(t)}dt,\ u\rw\IF, \ \lambda\rw\IF.
\EQN
Moreover,
\BQN\label{Thmeq22}
\Lambda_1(u)&\leq&\sum_{k=0}^{N(u)}\left(\mathcal{L}_u(I_{k}(u))+\mathcal{L}_u(I_{k+1}(u))
-\mathcal{L}_u(I_{k}(u)\cup I_{k+1}(u))\right)\nonumber\\
&\leq& \sum_{k=0}^{N(u)}\left(\pk{ \sup_{t\in I_{k}(u)} \LT\|\overline{\vk{X}}(t)*\vk{d}\RT\|_p^c>\mathcal{G}_{u,-\vn}(k)}+\pk{\sup_{t\in  I_{k+1}(u)}\LT\|\overline{\vk{X}}(t)*\vk{d}\RT\|_p^c>\mathcal{G}_{u,-\vn}(k)}\RT.\nonumber\\
& &\LT.-\pk{\sup_{t\in ((I_{k}(u)\cup I_{k+1}(u)))} \LT\|\overline{\vk{X}}(t)*\vk{d}\RT\|_p^c>\widehat{\mathcal{G}}_{u,-\vn}(k)}\right)\nonumber\\
&\leq &  \LT(2\mathcal{H}_\alpha[0,a^{1/\alpha}d^{-2/\alpha}\lambda]
-\mathcal{H}_\alpha[0,2a^{1/\alpha}d^{-2/\alpha}\lambda]\RT)
\sum_{k=0}^{N(u)}\pk{Z^c(0)>\widehat{\mathcal{G}}_{u,-\vn}(k)}
\nonumber\\
&\sim&
\frac{2\mathcal{H}_\alpha[0,a^{1/\alpha}d^{-2/\alpha}\lambda]
-\mathcal{H}_\alpha[0,2a^{1/\alpha}d^{-2/\alpha}\lambda]}{\lambda}
\int_0^{\IF}\exp\LT(-(1-\vn-\epsilon)f(t)\RT)dt\nonumber\\
&&\times u^{\frac{2}{\alpha^*}-\frac{2}{\beta^*}}\pk{Z^c(0)>u}\nonumber\\
&=&o\left( u^{\frac{2}{\alpha^*}-\frac{2}{\beta^*}}\pk{Z^c(0)>u}\right), \ u\rw\IF, \lambda\rw \IF, \vn\rw 0,\epsilon\rw 0,
\EQN
where $\widehat{\mathcal{G}}_{u,-\vn}(k)=\min (\mathcal{G}_{u,-\vn}(k),\mathcal{G}_{u,-\vn}(k+1))$.
By \nelem{in1}, we have
\BQN\label{Thmeq32}
\Lambda_2(u)
&\leq&\sum_{0\leq k,l\leq N(u), l\geq k+2}\pk{ \sup_{t\in I_{k}(u)}\LT\|\overline{\vk{X}}(t)*\vk{d}\RT\|_p^c>\mathcal{G}_{u,-\vn}(k), \sup_{t\in I_{l}(u)} \LT\|\overline{\vk{X}}(t)*\vk{d}\RT\|_p^c>\mathcal{G}_{u,-\vn}(l)}\nonumber\\
&\leq&\sum_{0\leq k\leq N(u)}\sum_{l=2}^{N(u)}\pk{ \sup_{t\in I_{k}(u)}\LT\|\overline{\vk{X}}(t)*\vk{d}\RT\|_p^c>\mathcal{G}_{u,-\vn}(k), \sup_{t\in I_{k+l}(u)}\LT\|\overline{\vk{X}}(t)*\vk{d}\RT\|_p^c>\mathcal{G}_{u,-\vn}(k)}\nonumber\\
&\leq&\mathbb{Q}_7\LT(\sum_{k=0}^{N(u)}\pk{Z^c(0)>\mathcal{G}_{u,-\vn}(k)}
\RT)\sum_{l=1}^{\IF}\exp\LT(-(l \lambda)^\alpha/8\RT)\nonumber\\
&\leq&\mathbb{Q}_8\pk{Z^c(0)>u} u^{\frac{2}{\alpha^*}-\frac{2}{\beta^*}}\lambda\sum_{l=1}^{\IF}\exp\LT(-(l \lambda)^\alpha/8\RT)\nonumber\\
&=& o\LT(\pk{Z^c(0)>u}u^{\frac{2}{\alpha^*}-\frac{2}{\beta^*}}\RT),\ u\rw\IF,\ \lambda\rw\IF.
\EQN
Combing (\ref{Th21})-(\ref{Thmeq32}) with (\ref{Thmeq12}), we obtain
\BQN\label{Th2re1}
\mathcal{L}_u([0,\delta(u)])\sim \pk{Z^c(0)>u}a^{1/\alpha}d^{-2/\alpha}\mathcal{H}_\alpha u^{\frac{2}{\alpha^*}-\frac{2}{\beta^*}}
\int_0^{\IF}e^{-f(t)}dt,\ u\rw\IF.
\EQN
{\bf Case 2:} $\beta^*=\alpha^*$.
We consider that for $u$ large enough,
\BQN
\mathcal{L}_u(I_{0}(u))\leq\mathcal{L}_u([0,\delta(u)])\leq
\sum_{k=0}^{N(u)}\mathcal{L}_u(I_{k}(u)).\label{Th22}
\EQN
Using \eqref{pp1} of Lemma \ref{im1} with $u$ replaced by $u^{1/c}$ and \eqref{hu2}, we have that
\BQN\label{Thmeq62}
\mathcal{L}_u(I_{0}(u))
&=& \pk{\sup_{t\in [0,\lambda]} \LT(Z^c(t u^{-2/\alpha^*})+g(t u^{-2/\alpha^*})\RT)>u}\nonumber\\
&\geq& \pk{\sup_{t\in [0,\lambda]} \frac{\LT\|\overline{\vk{X}}(tu^{-2/\alpha^*})*\vk{d}\RT\|_p^c}{1+\frac{w(1+\vn)}{u}\abs{tu^{-2/\alpha^*}}^\gamma+(1+\vn) bc\abs{tu^{-2/\alpha^*}}^\beta}>u}\nonumber\\
&\geq& \pk{\sup_{t\in [0,\lambda]} \frac{\LT\|\overline{\vk{X}}(tu^{-2/\alpha^*})*\vk{d}\RT\|_p}
{1+(1+\vn+\epsilon)u^{-2/c}d^{2}f(t)}>u^{1/c}}\nonumber\\
&\sim&\E{\sup_{t\in[0,\lambda]}\exp\LT(\sqrt{\frac{2a}{d^{2}}}B_{\alpha}(t)-\frac{a}{d^{2}}\abs{t}^{\alpha}-(1+\vn+\epsilon)f(t)\RT)}
\pk{Z^c(0)>u}\nonumber\\
&\sim&\mathcal{P}_{\alpha,\frac{a}{d^{2}}}^{f}[0,\IF)
\pk{Z^c(0)>u},\ u\rw\IF, \vn\rw 0, \epsilon\rw 0, \lambda\rw\IF.
\EQN
Similarly,
\BQN\label{Thmeq623}
\mathcal{L}_u(I_{0}(u))\leq \mathcal{P}_{\alpha,\frac{a}{d^{2}}}^{f}[0,\IF)
\pk{Z^c(0)>u},\ u\rw\IF, \lambda\rw\IF.
\EQN
Moreover, by \nelem{im1},
\BQN\label{Thmeq82}
{\sum_{k=1}^{N(u)}}\mathcal{L}_u(I_{k}(u))
&\leq& {\sum_{k=1}^{N(u)}} \pk{\sup_{t\in I_k(u)}\LT\|\overline{\vk{X}}(t)*\vk{d}\RT\|_p^c>\mathcal{G}_{u,-\vn}(k)}\nonumber\\
&\sim& \mathcal{H}_\alpha[0,a^{1/\alpha}d^{-2/\alpha}\lambda]
\sum_{k=1}^{N(u)}\pk{Z^c(0)>\mathcal{G}_{u,-\vn}(k)}\nonumber\\
&\leq&  \mathcal{H}_\alpha[0,a^{1/\alpha}d^{-2/\alpha}\lambda]
\pk{Z^c(0)>u}
\sum_{k=1}^{N(u)}\exp\LT(-(1-\vn-\epsilon)f\LT(k \lambda\RT)\RT)\nonumber\\
&\leq&
\mathcal{H}_\alpha[0,a^{1/\alpha}d^{-2/\alpha}\lambda]
\pk{Z^c(0)>u}
\sum_{k=1}^{\IF}\exp\LT(-\mathbb{Q}_9(k \lambda)^{\gamma\wedge\beta}\RT)\nonumber\\
&\sim&
\mathbb{Q}_{10}\pk{Z^c(0)>u}
\lambda\exp\LT(-\mathbb{Q}_{11} \lambda^{\gamma\wedge\beta}\RT)\nonumber\\
&=&
o\LT(\pk{Z^c(0)>u}\RT),\ u\rw\IF,\ \lambda\rw\IF.
\EQN
Inserting \eqref{Thmeq62}, \eqref{Thmeq623}, and \eqref{Thmeq82}  into \eqref{Th22}, we have
\BQN\label{Th2re2}
\mathcal{L}_u([0,\delta(u)])\sim\mathcal{P}_{\alpha,\frac{a}{d^{2}}}^{f}[0,\IF)
\pk{Z^c(0)>u}, u\rw\IF.
\EQN
{\bf Case 3:} $\beta^*<\alpha^*$. Obviously,
\BQN\label{Thmeq92}
\mathcal{L}_u([0,\delta(u)])\geq\pk{ Z^c(0)>u}.
\EQN
For any $\vn_2\in (0,1)$,  $[0,\delta(u)]\subseteq [0,u^{-2/\alpha^*}\vn_2]$ when $u$ large enough. By \nelem{im1} and the fact that $\sup_{t\in[0,\delta(u)]}g(t)\leq 0$, we obtain
\BQNY
\mathcal{L}_u([0,\delta(u)])
&\leq& \pk{\sup_{t\in [0,u^{-2/\alpha^*}\vn_2]}\LT\|\overline{\vk{X}}(t)*\vk{d}\RT\|_p^c>u}\nonumber\\
&\sim&\mathcal{H}_{\alpha}[0,a^{1/\alpha}d^{-2/\alpha}\vn_2]
\pk{ Z^c(0)>u}\nonumber\\
&\sim&\pk{ Z^c(0)>u},\ u\rw\IF, \vn_2\rw0.
\EQNY
Together with (\ref{Thmeq92}), we get
\BQN\label{Th2re3}
\mathcal{L}_u([0,\delta(u)])\sim \pk{ Z^c(0)>u},\ u\rw\IF.
\EQN
Consequently, we have the results according to \eqref{Main1}, \eqref{Th2re1}, \eqref{Th2re2} and \eqref{Th2re3}.\\
For $t_0\in(0,T)$ and $t_0=T$, we just need to replace $[0,\delta(u)]$ as $[t_0-\delta(u),t_0+\delta(u)]$ and $[T-\delta(u),T]$.
Thus we complete the proof.
\QED

\prooftheo{Thm0}
For any $\theta>0$ and $\lambda>0$, set $\alpha^*=\alpha c$
\BQNY
&&I_k(\theta)=[k\theta,(k+1)\theta],\quad a_k=a(k\theta), \quad k\in \N,\quad N(\theta)=\LT\lfloor\frac{T}{\theta}\RT\rfloor,\\
&&J^k_l(u)=\LT[k\theta+lu^{-2/\alpha^*}\lambda,k\theta+(l+1)u^{-2/\alpha^*}\lambda\RT],\quad M(u)=\LT\lfloor\frac{\theta u^{2/\alpha^*}}{\lambda}\RT\rfloor.
\EQNY
We have
\BQNY
\sum_{k=0}^{N(\theta)-1}\LT(\sum_{l=0}^{M(u)-1}\mathcal{K}_u(J^k_l(u))\RT)
-\sum_{i=1}^{4}\mathcal{A}_i(u)\leq\mathcal{K}_u([0,T])\leq\sum_{k=0}^{N(\theta)}
\mathcal{K}_u(I_k(\theta))
\leq\sum_{k=0}^{N(\theta)}\LT(\sum_{l=0}^{M(u)}\mathcal{K}_u(J^k_l(u))\RT),
\EQNY
where
\BQNY
&&\mathcal{A}_i(u)=\sum_{(k_1,l_1,k_2,l_2)\in \mathcal{L}_i}\mathcal{K}_u(J^{k_1}_{l_1},J^{k_2}_{l_2}),\ i=1,2,3,4,
\EQNY
with
\BQNY
&&\mathcal{L}_{1}=\LT\{0\leq k_1= k_2\leq N(\theta)-1,0\leq l_1+1=l_2\leq M(u)-1\RT\},\\
&&\mathcal{L}_{2}=\LT\{0\leq k_1+1= k_2\leq N(\theta)-1, l_1=M(u), l_2=0\RT\},\\
&&\mathcal{L}_3=\LT\{0\leq k_1+1<k_2\leq N(\theta)-1,0\leq l_1,l_2\leq M(u)-1\RT\},\\
&&\mathcal{L}_4=\LT\{0\leq k_1\leq k_2\leq N(\theta)-1,k_2-k_1\leq 1,0\leq l_1,l_2\leq M(u)-1\RT\}\setminus\LT(\mathcal{L}_1\cup \mathcal{L}_2\RT).
\EQNY
By \nelem{im1}
\BQNY
\sum_{k=0}^{N(\theta)}\LT( \sum_{l=0}^{M(u)}\mathcal{K}_u(J^k_l(u))\RT)
&=&\sum_{k=0}^{N(\theta)}\LT(\sum_{l=0}^{M(u)}\pk{\sup_{t\in [0,\lambda]}Z^c(k\theta+lu^{-2/\alpha^*}\lambda+u^{-2/\alpha^*} t)>u}\RT) \\
&\leq& \sum_{k=0}^{N(\theta)}\LT(\sum_{l=0}^{M(u)}(a_k+\vn_\theta)^{\frac{1}{\alpha}}d^{-2/\alpha}\mathcal{H}_\alpha \lambda\pk{Z^c(0)>u}\RT)\\
&\sim& \LT(\sum_{k=0}^{N(\theta)}(a_k+\vn_\theta\theta)^{\frac{1}{\alpha}}\RT)d^{-2/\alpha}\mathcal{H}_\alpha u^{2/\alpha^*} \pk{Z^c(0)>u}\\
&\sim& \int_{0}^T(a(t))^{1/\alpha}d tu^{2/\alpha^*}d^{-2/\alpha}\mathcal{H}_\alpha \pk{Z^c(0)>u}, \ u\rw\IF,\ \lambda\rw\IF,\ \theta\rw 0.
\EQNY
Similarly,
\BQNY
\sum_{k=0}^{N(\theta)-1}\LT(\sum_{l=0}^{M(u)-1}\mathcal{K}_u(J^k_l(u))\RT)
\geq \int_{0}^T(a(t))^{1/\alpha}d tu^{2/\alpha^*}d^{-2/\alpha}\mathcal{H}_\alpha \pk{Z^c(0)>u}, \ u\rw\IF,\ \lambda\rw\IF,\ \theta\rw 0.
\EQNY
Further, by \nelem{im1}
\BQNY
\mathcal{A}_1(u)
&=&\sum_{k=0}^{N(\theta)-1}\LT(\sum_{l=0}^{M(u)-1}\left(\mathcal{K}_u(J^k_l(u))
+\mathcal{K}_u(J^k_{l+1}(u))-\mathcal{K}_u(J^k_l(u)\cup J^k_{l+1}(u))\right)\RT)\nonumber\\
&\sim &  \sum_{k=0}^{N(\theta)-1}\LT(\LT(\mathcal{H}_\alpha[0,(a_k+\vn_\theta)^{\frac{1}{\alpha}}d^{-1/\alpha}\lambda]
+\mathcal{H}_\alpha[0,(a_k+\vn_\theta)^{\frac{1}{\alpha}}d^{-1/\alpha}\lambda]
-\mathcal{H}_\alpha[0,2(a_k-\vn_\theta)^{\frac{1}{\alpha}}d^{-1/\alpha}\lambda]\RT)\RT.\\
&&\LT.\times\sum_{l=0}^{M(u)-1}\pk{Z^c(0)>u}\RT)\nonumber\\
&\leq&\mathbb{Q}_1\LT(\sum_{k=0}^{N(\theta)-1}\LT((a_k+\vn_\theta)^{\frac{1}{\alpha}}-(a_k
-\vn_\theta)^{\frac{1}{\alpha}}\RT)\theta\RT)u^{2/\alpha^*}\pk{Z^c(t)>u}\\
&=&o\left(u^{2/\alpha^*}\pk{Z^c(t)>u}\right), \ u\rw\IF,\ \ \lambda\rw\IF,\theta\rw0.
\EQNY
Similarly, by \nelem{im1}
\BQNY
\mathcal{A}_2(u)
&=&\sum_{k=0}^{N(\theta)-1}\mathcal{K}_u(J^{k}_{M(u)-1}(u), J^{k+1}_{0}(u))\\
&\leq&\sum_{k=0}^{N(\theta)-1}\pk{\sup_{t\in [0,2\lambda] }Z^c((k+1)\theta-u^{-2/\alpha^*}t)>u,\sup_{t\in [0,2\lambda]}Z^c((k+1)\theta+u^{-2/\alpha^*}t)>u}\\
&=&\sum_{k=0}^{N(\theta)-1}\LT(\pk{\sup_{t\in [0,2\lambda] }Z^c((k+1)\theta-u^{-2/\alpha^*}t)>u}+\pk{\sup_{t\in [0,2\lambda]}Z^c((k+1)\theta+u^{-2/\alpha^*}t)>u}\RT.\\
&&\LT.-\pk{\sup_{t\in [-2\lambda,2\lambda] }Z^c((k+1)\theta-u^{-2/\alpha^*}t)>u}\RT)\nonumber\\
&\sim &  \sum_{k=0}^{N(\theta)-1}\LT(\LT(2\mathcal{H}_\alpha[0,2(a_{k+1}+\vn_\theta)^{\frac{1}{\alpha}}d^{-1/\alpha}\lambda]
-\mathcal{H}_\alpha[-2(a_k-\vn_\theta)^{\frac{1}{\alpha}}d^{-1/\alpha}\lambda,2(a_k-\vn_\theta)^{\frac{1}{\alpha}}d^{-1/\alpha}\lambda]\RT)\RT.\\
&&\LT.\times\sum_{l=0}^{M(u)-1}\pk{Z^c(0)>u}\RT)\nonumber\\
&\leq&\mathbb{Q}_2\LT(\sum_{k=0}^{N(\theta)-1}\LT((a_k+\vn_\theta)^{\frac{1}{\alpha}}-(a_k
-\vn_\theta)^{\frac{1}{\alpha}}\RT)\theta\RT)u^{2/\alpha^*}
\pk{Z^c(0)>u}\\
&=&o\left(u^{2/\alpha^*}\pk{Z^c(0)>u}\right), \ u\rw\IF,\ \ \lambda\rw\IF,\theta\rw0.
\EQNY

For any $\theta>0$
\BQNY
\E{X_i(t)X_i(s)}=r(s,t)\leq 1-\delta(\theta)
\EQNY
for $(s,t)\in J^{k_1}_{l_1}(u)\times J^{k_2}_{l_2}(u),(j_1,k_1,j_2,k_2)\in\mathcal{L}_3$ where $\delta(\theta)>0$ is related to $\theta$. Then by \nelem{tt}
\BQNY
\mathcal{A}_3(u)
&\leq&N(\theta)M(u)2\Psi\LT(\frac{2\LT(u\RT)^{\frac{1}{c}}-\mathbb{Q}_3}{d\sqrt{4-\delta(\theta)}}\RT)\\
&\leq&\frac{T}{\lambda}u^{2/\alpha^*}2\Psi\LT(\frac{2u^{\frac{1}{c}}-\mathbb{Q}_3}{d\sqrt{4-\delta(\theta)}}\RT)\\
&=&o\LT(u^{2/\alpha^*}\pk{Z^c(0)>u}\RT), \ u\rw\IF, \lambda\rw\IF, \theta\rw 0.
\EQNY
where $\mathbb{Q}_3$ is a large constant.
Finally by \nelem{in1} for $u$ large enough and $\theta$ small enough
\BQNY
\mathcal{A}_4(u)
&\leq &\sum_{k=0}^{N(\theta)-1}\LT(\sum_{l=0}^{2M(u)}\sum_{i=2}^{2M(u)}
\mathcal{K}_u(J^k_l(u),J^{k}_{l+i}(u))\RT)\\
&\leq&\sum_{k=0}^{N(\theta)-1} \sum_{l=0}^{2M(u)}\pk{Z^c(0)>u}
\LT(\sum_{i=1}^{\IF}\mathbb{Q}_4\exp\LT(-\frac{\mathbb{Q}_5}{8}\abs{i\lambda}^\alpha\RT)\RT)\\
&\leq&\mathbb{Q}_6\frac{T}{\lambda}u^{-2/\alpha^*}\pk{Z^c(0)>u}\LT(\sum_{i=1}^{\IF}\exp\LT(-\frac{\mathbb{Q}_5}{8}\abs{i\lambda}^\alpha\RT)\RT)\\
&=&o\LT(u^{2/\alpha^*}\pk{Z^c(0)>u}\RT), \ u\rw\IF, \lambda\rw\IF, \theta\rw 0.
\EQNY
Thus the claim follows.
\QED

\prooftheo{Thm1}
Through this proof, denote $\mathcal{L}_u(\Delta_1)$ and $\mathcal{L}_u(\Delta_1,\Delta_2)$ the same as in the proof of \netheo{Thm2}.\\
When $c\in(0,2)$, in the proof of \netheo{Thm2}, if we take $\beta^*=\frac{2\gamma c}{2-c}$ and $f(t)=\frac{w\abs{t}^\gamma}{cd^2}$, then all argumentations still hold and the results follow.

When $c=2$, for any constant $\theta>0$, we define
$$I_k=[k\theta, (k+1)\theta],\ k\in\N, \  N(\theta)= \LT\lfloor\frac{T}{\theta}\RT\rfloor,$$
and $$M_1(k)=\sup_{t\in I_k}g(t),\ \ \  M_2(k)=\inf_{t\in I_k}g(t). $$
Then
\BQNY
\mathcal{L}_u([0,T])\geq\sum_{k=0}^{N(\theta)-1}\mathcal{L}_u(I_k)-\sum_{j=1}^2\Lambda_j,
\EQNY
where
\BQNY
\Lambda_1=\sum_{k=0}^{N(\theta)}\mathcal{L}_u(I_k, I_{k+1}),\ \
\Lambda_2=\underset{j>k+1}{\sum_{k=0}^{N(\theta)}}\mathcal{L}_u(I_k, I_{j}),
\EQNY
and by \netheo{Thm0}
\BQNY
\mathcal{L}_u([0,T])
&\leq&\sum_{k=0}^{N(\theta)}\mathcal{L}_u(I_k) \\
&\leq&\sum_{k=0}^{N(\theta)}\pk{\sup_{t\in I_k}Z^c(t)>u-M_1(k)}\\
&\sim& \sum_{k=0}^{N(\theta)}(a(k\theta))^{\frac{1}{\alpha}}\LT(u-M_1(k)\RT)^{\frac{1}{\alpha }}d^{-2/\alpha}\mathcal{H}_\alpha \theta\pk{Z^c(0)>u-M_1(k)}\\
&\sim& u^{\frac{1}{\alpha}}d^{-2/\alpha}\mathcal{H}_\alpha \pk{Z^c(0)>u}\theta\sum_{k=0}^{N(\theta)}(a(k\theta))^{\frac{1}{\alpha}}e^{\frac{M_1(k)}{2d^2}}\\
&\sim& u^{\frac{1}{\alpha}}d^{-2/\alpha}\mathcal{H}_\alpha  \pk{Z^c(0)>u}\int_{0}^T(a(t))^{1/\alpha}e^{\frac{g(t)}{2d^2}}d t, \ u\rw\IF,\ \theta\rw 0.
\EQNY
Similarly,
\BQNY
\sum_{k=0}^{N(\theta)-1}\mathcal{L}_u(I_k)
&\geq&\sum_{k=0}^{N(\theta)-1}\pk{\sup_{t\in I_k}Z^c(t)>u-M_2(k)}\\
&\sim& u^{\frac{1}{\alpha}}d^{-2/\alpha}\mathcal{H}_\alpha  \pk{Z^c(0)>u}\int_{0}^T(a(t))^{1/\alpha}e^{\frac{g(t)}{2d^2}}d t, \ u\rw\IF,\ \theta\rw 0.
\EQNY
Further, we have
\BQNY
\Lambda_1&\leq&\sum_{k=0}^{N(\theta)}\left(\mathcal{L}_u(I_k)+\mathcal{L}_u(I_{k+1})
-\mathcal{L}_u(I_{k}\cup I_{k+1})\right)\nonumber\\
&\leq&\sum_{k=0}^{N(\theta)}\left(\pk{ \sup_{t\in I_{k}}  Z^c(t)>u-\widetilde{M}_1(k)}+\pk{\sup_{t\in I_{k+1}} Z^c(t)>u-\widetilde{M}_1(k)}\RT.\nonumber\\
&&\LT.-\pk{\sup_{t\in I_{k}\cup I_{k+1}}  Z^c(t)>u-\widetilde{M}_1(k)}\right)\nonumber\\
&\sim &  \sum_{k=0}^{N(\theta)}\LT((a(k\theta))^{1/\alpha}+(a((k+1)\theta))^{1/\alpha}
-2(a(k\theta))^{1/\alpha}\RT)\theta
u^{\frac{1}{\alpha }}d^{-2/\alpha}\mathcal{H}_\alpha e^{\frac{\widetilde{M}_1(k)}{2d^2}}\pk{Z^c(0)>u}
\nonumber\\
&=&o\left(u^{1/\alpha}\pk{Z^c(0)>u}\right), \ u\rw\IF,\ \theta\rw 0,
\EQNY
where $\widetilde{M}_1(k)=\max(M_1(k),M_1(k+1))$.\\
Then for $g_m=\sup_{t\in[0,T]}g(t)$ by \nelem{tt}
\BQNY
\Lambda_2
&\leq& \underset{j>k+1}{\sum_{k=0}^{N(\theta)}}\pk{\sup_{t\in I_k}Z^c(t)>u-g_m,\sup_{t\in I_{j}}Z^c(t)>u-g_m}\\
&\leq&\underset{j>k+1}{\sum_{k=0}^{N(\theta)}}2\Psi\LT(\frac{2{(u-g_m)^{1/c}}-\mathbb{Q}_1}{d\sqrt{4-\delta(\theta)}}\RT)\\
&=&o\LT(\pk{Z^c(0)>u}\RT), \ u\rw\IF, \theta\rw0.
\EQNY
Thus, we have
\BQNY
\mathcal{L}_u([0,T])\sim \int_{0}^T(a(t))^{1/\alpha}e^{\frac{g(t)}{2d^2}}d t \mathcal{H}_\alpha d^{-2/\alpha}  u^{\frac{1}{\alpha}}\pk{Z^c(0)>u}, \ u\rw\IF.
\EQNY
When $c\in(2,\IF)$, set $M_1=\inf_{t\in[0,T]}g(t)$ and $M_2=\sup_{t\in[0,T]}g(t)$. Since $g(t)$ is a continuous function, we have
$-\IF<M_1\leq M_2<\IF$.
Further, since when $c\in(2,\IF)$, $$\pk{Z^c(0)>u+\mathbb{Q}_2}\sim\pk{Z^c(0)>u}$$
holds for any $\mathbb{Q}_2>0$.
Hence, by  \netheo{Thm0}
\BQNY
\mathcal{L}_u([0,T])&\geq& \pk{\sup_{t\in[0,T]}Z^c(t)>u-M_1}\\
&\sim& \int_{0}^T(a(t))^{\frac{1}{\alpha}}dt d^{-\frac{2}{\alpha}} \mathcal{H}_\alpha u^{\frac{2}{\alpha c}}\pk{Z^c(0)>u-M_1}\\
&\sim& \int_{0}^T(a(t))^{\frac{1}{\alpha}}dt d^{-\frac{2}{\alpha}} \mathcal{H}_\alpha u^{\frac{2}{\alpha c}}\pk{Z^c(0)>u}, u\rw\IF,
\EQNY
and
\BQNY
\mathcal{L}_u([0,T])\leq \pk{\sup_{t\in[0,T]}Z^c(t)>u-M_2}
\sim \int_{0}^T(a(t))^{\frac{1}{\alpha}}dt d^{-\frac{2}{\alpha}} \mathcal{H}_\alpha u^{\frac{2}{\alpha c}}\pk{Z^c(0)>u}, u\rw\IF.
\EQNY
The result follows.
\QED

\section{Appendix}
In this section, we give several lemmas which are used in the proofs of the theorems.
\BEL\label{tt}
let $\vk{X}(t)=(X_1(t)\ldots,X_n(t)),\ t\in[0,T],n\ge 1$ be an centered ${\R}^n$-valued  vector process with independent marginals, which have continuous samples, unit variances and correlation functions satisfying assumption (v).
 Then for $ 0<t_1<t_2<t_3<\IF$ and $u$ large enough
\BQNY
\pk{\sup_{t\in[0,t_1]} Z^c(t)>u,\sup_{t\in[t_2,t_3]} Z^c(t)>u}\leq
2\Psi\LT(\frac{2u^{1/c}-D}{d\sqrt{4-\delta}}\RT),
\EQNY
where $D, \delta$ are some constant.
\EEL
\prooflem{tt}
By assumption (v) and the continuity of $r(t)$, for some $\delta>0$ we have
\BQNY
\E{X_i(t)X_i(s)}=r(s,t)\leq1-\frac{\delta}{2}, i=1,2,\ldots,n,
\EQNY
holds for any $(s,t)\in[0,t_1]\times[t_2,t_3]$.
Set $\widetilde{Y}(t,\vk{v},s, \vk{w})=\sum_{i=1}^n X_i(t)d_iv_i+\sum_{i=1}^n X_i(s)d_iw_i$ where $\vk{v},\vk{w}\in\mathcal{S}_{q} $  with $\mathcal{S}_{q}=\{\vk{v}\in\R^n:||\vk{v}||_q=1\}$. %,\mathbb{S}^2_{d-1}$ two unit spheres.
Since  $\widetilde{Y}(t,\vk{v},s, \vk{w})$ is a center Gaussian fields, we have further
\BQNY
\Var\LT(\widetilde{Y}(t,\vk{v},s, \vk{w})\RT)&=&\sum_{i=1}^n(v_i^2+w_i^2+2r(s,t)v_iw_i)d_i^2\\
&\leq&2d^2+2r(s,t)\sum_{i=1}^n(v_i^2+w_i^2)d_i^2\\
&=&2d^2+2d^2r(s,t)\\
&\leq&d^2(4-\delta),
\EQNY
for any $(t,\vk{v},s,\vk{w})\in[0,t_1]\times\mathcal{S}_{q}\times[t_2,t_3]\times\mathcal{S}_{q}$.
By Borell inequality,
\BQNY
\pk{\sup_{t\in[0,t_1]} Z^c(t)>u,\sup_{t\in[t_2,t_3]} Z^c(t)>u}
&=&\pk{\sup_{t\in[0,t_1]} Z(t)>u^{1/c},\sup_{t\in[t_2,t_3]} Z(t)>u^{1/c}}\\
&\leq&\pk{\sup_{(t,\vk{v},s,\vk{w})\in[0,t_1]\times\mathcal{S}_{q}\times[t_2,t_3]\times\mathcal{S}_{q}}\widetilde{Y}(t,\vk{v},s, \vk{w})>2u^{1/c}}\\
&\leq& 2\Psi\LT(\frac{2u^{1/c}-D}{d\sqrt{4-\delta}}\RT),
\EQNY
where $D$ is some constant such that
$$\pk{\sup_{(t,\vk{v},s,\vk{w})\in[0,t_1]\times\mathcal{S}_{q}\times[t_2,t_3]\times\mathcal{S}_{q}}\widetilde{Y}(t,\vk{v},s, \vk{w})>D}\leq \frac{1}{2},$$
 hence the claim follows.
\QED

\BEL\label{im1}
let $\vk{X}(t)=(X_1(t)\ldots,X_n(t)),\ t\in\R,n\ge 1$ be an centered ${\R}^n$-valued  vector process with independent marginals, which have continuous samples, unit variances and correlation functions satisfying assumption (iv).
Set $ a:=a(t_0), \ t_0\in\R,$ and $K_u$ a family of index sets and $u_k$ satisfying that
\BQN\label{uk}
\lim_{u\rw\IF}\sup_{k\in K_u}\LT|\frac{u_k}{u}-1\RT|=0.
\EQN
If $f(t)$ is a nonnegative  continuous function with $f(0)=0,f(t)>0,t\neq 0$ and $\vk{d}$ is the same as in \eqref{vd},
then we have that for some constants $S_1, S_2\geq 0$ and $\max(S_1,S_2)>0$
\BQN\label{pp1}
{\pk{\sup_{t\in [-S_1,S_2]}\frac{Z(u^{-2/\alpha}t+t_0 )}{1+u^{-2}f(t)}>u}}
\sim\mathcal{P}_{\alpha,ad^{-2}}^{\frac{1}{d^{2}}f(t)}
[-S_1,S_2]\pk{Z(t_0)>u},\ u\rw\IF,
\EQN
and
\begin{align*}\label{pp2}
\mathcal{P}_{\alpha,ad^{-2}}^{\frac{1}{d^{2}}f(t)}
[-S_1,S_2]=\E{\exp\LT(\sup_{t\in[-S_1,S_2]}
\sqrt{\frac{2a}{d^{2}}}B_{\alpha}(t)-\frac{a}{d^{2}}\abs{t}^\alpha-\frac{1}{d^{2}}f(t)\RT)}.
\end{align*}
If $\lim_{u\rw\IF}\sup_{k\in K_u}\abs{k u^{-2/\alpha}}\leq \theta$ for some small enough $\theta\geq0$,
we have for some constant $S>0$
\BQN\label{P2p2}
\mathcal{H}_\alpha[-(a-\vn_{\theta})^{1/\alpha}d^{-2/\alpha}S_1,
(a-\vn_{\theta})^{1/\alpha}d^{-2/\alpha}S_2]
&\leq&\lim_{u\rw\IF}\forall_{k\in K_u}\frac{\pk{\sup_{t\in [-S_1,S_2]}Z(u^{-\frac{2}{\alpha}}(t+kS)+t_0 )>u_k}}
{\pk{Z(t_0)>u_k}}\nonumber\\
&\leq& \mathcal{H}_\alpha[-(a+\vn_{\theta})^{1/\alpha}d^{-2/\alpha}S_1,
(a+\vn_{\theta})^{1/\alpha}d^{-2/\alpha}S_2],
\EQN
where $\vn_\theta\rw 0$, as $\theta\rw 0$.\\
Specially, if $\theta=0$, we have
\BQN\label{P2p}
\lim_{u\rw\IF}\sup_{k\in K_u}\LT|\frac{\pk{\sup_{t\in [-S_1,S_2]}Z(u^{-\frac{2}{\alpha}}(t+kS)+t_0 ) >u_k}}
{\pk{Z(t_0)>u_k}}-\mathcal{H}_\alpha[-a^{\frac{1}{\alpha}}d^{-\frac{2}{\alpha}}S_1,a^{\frac{1}{\alpha}}d^{-\frac{2}{\alpha}}S_2]\RT|=0,
\EQN
and
\begin{align*}
\mathcal{H}_\alpha[-a^{1/\alpha}d^{-2/\alpha}S_1,a^{1/\alpha}d^{-2/\alpha}S_2]=\E{\exp\LT(\sup_{t\in[-a^{1/\alpha}d^{-2/\alpha}S_1,a^{1/\alpha}d^{-2/\alpha}S_2]}\sqrt{2}B_{\alpha}(t)-|t|^\alpha\RT)}.
\end{align*}
\EEL
\prooflem{im1}
{\bf Step 1:}
First we give the proof of \eqref{pp1}.
When $p=1$, set $\mathcal{W}=\{\vk{w}=(w_1,\cdots,w_n): w_i=\pm 1,\ i=1,\cdots, n\}$.
Then we have
\BQNY
&&\pk{\sup_{t\in [-S_1,S_2]}\frac{Z(u^{-2/\alpha}t+t_0 ) }{1+u^{-2}f(t)}>u}\\
&&=\pk{\sup_{t\in [-S_1,S_2]}\frac{\sum_{i=1}^n\abs{d_i X_i(u^{-2/\alpha}t+t_0)}  }{1+u^{-2}f(t)}>u}\\
&&=\sum_{\vk{w}\in\mathcal{W}}\pk{\sup_{t\in [-S_1,S_2]}\frac{\sum_{i=1}^nw_id_iX_i(u^{-2/\alpha}t+t_0)  }{1+u^{-2}f(t)}>u}\\
&&\quad-\underset{\vk{w}\neq \vk{w}'}{\sum_{\vk{w},\vk{w}'\in\mathcal{W}}}
\pk{\sup_{t\in [-S_1,S_2]}\frac{\sum_{i=1}^nw_id_iX_i(u^{-2/\alpha}t+t_0)  }{1+u^{-2}f(t)}>u,\sup_{s\in [-S_1,S_2]}\frac{\sum_{i=1}^nw'_id_iX_i(u^{-2/\alpha}s+t_0)  }{1+u^{-2}f(s)}>u}\\
&&=2^n\pk{\sup_{t\in [-S_1,S_2]}\frac{\sum_{i=1}^nd_iX_i(u^{-2/\alpha}t+t_0)  }{1+u^{-2}f(t)}>u}\\
&&\quad-\underset{\vk{w}\neq \vk{w}'}{\sum_{\vk{w},\vk{w}'\in\mathcal{W}}}
\pk{\sup_{t\in [-S_1,S_2]}\frac{\sum_{i=1}^nw_id_iX_i(u^{-2/\alpha}t+t_0)  }{1+u^{-2}f(t)}>u,\sup_{s\in [-S_1,S_2]}\frac{\sum_{i=1}^nw'_id_iX_i(u^{-2/\alpha}s+t_0)  }{1+u^{-2}f(s)}>u}.
\EQNY
By \cite{Threshold2016} [Lemma 4.1], we have
\BQNY
2^n\pk{\sup_{t\in [-S_1,S_2]}\frac{\sum_{i=1}^nd_iX_i(u^{-2/\alpha}t+t_0)  }{1+u^{-2}f(t)}>u}&=&2^n\pk{\sup_{t\in [-S_1,S_2]}\frac{\sum_{i=1}^nd_iX_i(u^{-2/\alpha}t+t_0)  }{1+u^{-2}f(t)}>u}\\
&=&2^n\pk{\sup_{t\in [-S_1,S_2]}\frac{\LT(\sum_{i=1}^nd^2_i\RT)^{1/2}X_1(u^{-2/\alpha}t+t_0)  }{1+u^{-2}f(t)}>u}\\
&\sim&2^{n}\mathcal{P}_{\alpha,ad^{-2}}^{\frac{1}{d^{2}}f(t)}
[-S_1,S_2]\Psi\LT(\frac{u}{d}\RT)\\
&\sim&\mathcal{P}_{\alpha,ad^{-2}}^{\frac{1}{d^{2}}f(t)}
[-S_1,S_2]\pk{Z(t_0)>u},\ u\rw\IF.
\EQNY
Since for any $\vk{w}\neq \vk{w}'$
\BQNY
V^2_1:&=&\E{\LT(\sum_{i=1}^nw_i d_iX_i(u^{-2/\alpha}t+t_0)+\sum_{i=1}^nw'_i d_iX_i(u^{-2/\alpha}s+t_0)\RT)^2}\\
&=&2\sum_{i=1}^nd_i^2
+2\sum_{i=1}^nw_iw_i'd_i^2r(u^{-2/\alpha}t+t_0,u^{-2/\alpha}s+t_0)\\
&<& 4\sum_{i=1}^nd_i^2=4d^2,
\EQNY
then by Borell inequality, we have
\BQNY
&&\pk{\sup_{t\in [-S_1,S_2]}\frac{\sum_{i=1}^nw_id_iX_i(u^{-2/\alpha}t+t_0)  }{1+u^{-2}f(t)}>u,\sup_{s\in [-S_1,S_2]}\frac{\sum_{i=1}^nw'_id_iX_i(u^{-2/\alpha}s+t_0)  }{1+u^{-2}f(s)}>u}\\
&&\leq\pk{\sup_{(t,s)\in [-S_1,S_2]\times[S_1,S_2]}\LT(\sum_{i=1}^nw_id_iX_i(u^{-2/\alpha}t+t_0)+  \sum_{i=1}^nw'_id_iX_i(u^{-2/\alpha}s+t_0)\RT)>2u}\\
&&\leq\exp\LT(-\frac{(2u-\mathbb{Q})^2}{2V^2_1}\RT)\\
&&=o\LT(\pk{Z(t_0)>u}\RT),\ u\rw\IF.
\EQNY
Then \eqref{pp1} with $p=1$ is follow.\\
When $p\in(1,\IF]$, set $Y(t,\vk{v})=\sum_{i=1}^n d_iv_iX_{i}(t),(t,\vk{v})\in\R\times\mathcal{S}_q$ which is a centered Gaussian field.\\
Then we have
\BQNY
\pk{\sup_{t\in [-S_1,S_2]}\frac{Z(u^{-2/\alpha}t+t_0 ) }{1+u^{-2}f(t)}>u}
=\pk{\sup_{(t,\vk{v})\in [-S_1,S_2]\times\mathcal{S}_q}\frac{Y(u^{-2/\alpha}t+t_0 ,\vk{v})}{1+u^{-2}f(t)}>u}.
\EQNY
Set $\mathcal{S}_q^{\delta}=\LT\{\vk{v}\in\mathcal{S}_q:d^2-\sum_{i=1}^nd_i^2v_i^2\leq \delta\RT\},\delta>0$.
Next we prove that as $u\rw\IF$
\BQNY
\pk{\sup_{(t,\vk{v})\in [-S_1,S_2]\times\mathcal{S}_q}\frac{Y(u^{-2/\alpha}t+t_0 ,\vk{v})}{1+u^{-2}f(t)}>u}\sim\pk{\sup_{(t,\vk{v})\in [-S_1,S_2]\times\mathcal{S}^\delta_q}\frac{Y(u^{-2/\alpha}t+t_0 ,\vk{v})}{1+u^{-2}f(t)}>u}.
\EQNY
Since
\BQNY
\pk{\sup_{(t,\vk{v})\in [-S_1,S_2]\times\mathcal{S}^\delta_q}\frac{Y(u^{-2/\alpha}t+t_0 ,\vk{v})}{1+u^{-2}f(t)}>u}\geq\pk{\sup_{\vk{v}\in\mathcal{S}^\delta_q}Y(t_0 ,\vk{v})>u}=\pk{Z(t_0)>u},
\EQNY
we just need to show as $u\rw\IF$
\BQNY
\pk{\sup_{(t,\vk{v})\in [-S_1,S_2]\times(\mathcal{S}_q\setminus\mathcal{S}^\delta_q)}\frac{Y(u^{-2/\alpha}t+t_0 ,\vk{v})}{1+u^{-2}f(t)}>u}=o\LT(\pk{Z(t_0)>u}\RT).
\EQNY
In fact, since
\BQNY
\sup_{(t,\vk{v})\in[-S_1,S_2]\times(\mathcal{S}_q\setminus\mathcal{S}^\delta_q)}
\Var(Y(u^{-2/\alpha}t+t_0,\vk{v}))=\sup_{\vk{v}\in(\mathcal{S}_q\setminus\mathcal{S}^\delta_q)}
\LT(\sum_{i=1}^n d_i^2v_i^2\RT)\leq d^2-\delta,
\EQNY
by Borell inequality, we have
\BQNY
\pk{\sup_{(t,\vk{v})\in [-S_1,S_2]\times(\mathcal{S}_q\setminus\mathcal{S}^\delta_q)}\frac{Y(u^{-2/\alpha}t+t_0 ,\vk{v})}{1+u^{-2}f(t)}>u}
&\leq& \pk{\sup_{(t,\vk{v})\in [-S_1,S_2]\times(\mathcal{S}_q\setminus\mathcal{S}^\delta_q)}{Y(u^{-2/\alpha}t+t_0 ,\vk{v})}>u}\\
&\leq& \exp\LT(-\frac{\LT(u-\mathbb{Q}_1\RT)^2}{2(d^2-\delta)}\RT)\\
&=&o\LT(\pk{Z(t_0)>u}\RT),\ u\rw\IF,
\EQNY
where $\mathbb{Q}_1:=\E{\sup_{(t,\vk{v})\in[-S_1,S_2]\times(\mathcal{S}_q\setminus\mathcal{S}^\delta_q)}
Y(u^{-2/\alpha}t+t_0 ,\vk{v})}<\IF$.\\
When $p\in(1,2)\cup(2,\IF]$, by \nelem{dd}, we know $\sigma_1^2(t,\vk{v}):=\Var\LT(\frac{Y(u^{-2/\alpha}t+t_0 ,\vk{v})}{1+u^{-2}f(t)}\RT)$ attains the maximum over $[-S_1,S_2]\times\mathcal{S}_q$ at several discrete points, so we can choose $\delta$ small enough such that
$\mathcal{D}^i_\delta=[-S_1,S_2]\times\mathcal{S}^\delta_q(i) $ with $\mathcal{S}^\delta_q(i)$ the union of non-overlapping compact neighborhoods of $ \vk{v}_+^i, \vk{v}_-^i$ or $\vk{z}$ in \nelem{dd}.
Then as mentioned in \cite{Pit96} or \cite{Fatalov1992}[Lemma 2.1]
\BQN\label{sum1}
\pk{\sup_{(t,\vk{v})\in [-S_1,S_2]\times\mathcal{S}_q}\frac{Y(u^{-2/\alpha}t+t_0 ,\vk{v})}{1+u^{-2}f(t)}>u}\sim
\sum_{i=1}^M\pk{\sup_{(t,\vk{v})\in \mathcal{D}^i_\delta}\frac{Y(u^{-2/\alpha}t+t_0 ,\vk{v})}{1+u^{-2}f(t)}>u},\ u\rw\IF,
\EQN
where $M$ is the number of the maximum point of $\sigma_1^2(t,\vk{v})$.\\
Case 1) $p\in(1,2)$ and $M=2^n$. It is enough to find the asymptotics of single term in \eqref{sum1}, for instance, for a point $(0,\vk{z}),\ z_i=(d_i/d)^{2/q-2}$. In a neighborhood $\mathcal{S}^\delta_q(1)$ of $\vk{z}$, we have
$$v_n=\LT(1-\sum_{i=1}^{n-1}v_i^q\RT)^{1/q},$$
hence the fields $\frac{Y(u^{-2/\alpha}t+t_0 ,\vk{v})}{1+u^{-2}f(t)}$ can be represented as
\BQNY
Y_1(u^{-2/\alpha}t+t_0, \widetilde{\vk{v}})
=\sum_{i=1}^{n-1}v_id_i \frac{X_i(u^{-2/\alpha}t+t_0)}{1+u^{-2}f(t)}
+\LT(1-\sum_{i=1}^{n-1}v_i^q\RT)^{1/q}d_n\frac{X_n(u^{-2/\alpha}t+t_0)}{1+u^{-2}f(t)},
 \widetilde{\vk{v}}=(v_1,\cdots,v_{n-1}),
\EQNY
which is defined in $[-S_1,S_2]\times\widetilde{\mathcal{S}}^\delta_q(1)$ where
$$\widetilde{\mathcal{S}}^\delta_q(1)
=\LT\{\widetilde{\vk{v}}:\LT(v_1,\cdots,v_{n-1},\LT(1-\sum_{i=1}^{n-1}v_i^q\RT)^{1/q}\RT)\in\mathcal{S}^\delta_q(1)\RT\},$$
is a small neighborhood of $\widetilde{\vk{z}}=\LT(z_1,\cdots,z_{n-1}\RT)$.
On $[-S_1,S_2]\times\widetilde{\mathcal{S}}^\delta_q(1)$, the variance
\BQNY
\sigma_1^2(t,\widetilde{\vk{v}}):=\frac{1}{(1+u^{-2}f(t))^2}\sigma_1^2(\widetilde{\vk{v}})
:=\frac{1}{(1+u^{-2}f(t))^2}\LT[\sum_{i=1}^{n-1}d_i^2v_i^2+d_n^2\LT(1-\sum_{i=1}^{n-1}v_i^q\RT)^{2/q}\RT]
\EQNY
of $Y_1(u^{-2/\alpha}t+t_0, \widetilde{\vk{v}})$ attains its maximum $d^2$ at $(0,\widetilde{\vk{z}})$ where $\widetilde{\vk{z}}$ is a interior point of a set $\widetilde{\mathcal{S}}^\delta_q(1)$. We can write the following Taylor expansion for $\sigma_1(t,\widetilde{\vk{v}})$
\BQNY
\sigma_1(t,\widetilde{\vk{v}})=\frac{d}{1+u^{-2}f(t)}
-\frac{q-2}{2d}(\widetilde{\vk{v}}-\widetilde{\vk{z}})
\Lambda(\widetilde{\vk{v}}-\widetilde{\vk{z}})^T+o(|\widetilde{\vk{v}}-\widetilde{\vk{z}}|^2), \widetilde{\vk{v}}\rw\widetilde{\vk{z}},\ u\rw\IF,
\EQNY
where $\Lambda=(\lambda_{i,j})_{i,j=1,\cdots,n-1}$ is a non-negative define matrix with elements
\BQNY
\lambda_{i,j}=-(2(q-2))^{-1}\frac{\partial^2}{\partial v_i\partial v_j}\LT[\sum_{i=1}^{n-1}d_i^2v_i^2+d_n^2\LT(1-\sum_{i=1}^{n-1}v_i^q\RT)^{2/q}\RT]|_{\widetilde{\vk{v}}=\widetilde{\vk{z}}},
i,j=1,\cdots,n-1.
\EQNY
We have the following expansion for the correlation function $r_1(t,\widetilde{\vk{v}},s,\widetilde{\vk{v}}')$ of $Y_1(u^{-2/\alpha}t+t_0, \widetilde{\vk{v}})$
\BQNY
r_1(t,\widetilde{\vk{v}},s,\widetilde{\vk{v}}')=1- u^{-2}a(t-s)^\alpha
-\frac{1}{2d}(\widetilde{\vk{v}}-\widetilde{\vk{v}}')
\Lambda(\widetilde{\vk{v}}-\widetilde{\vk{v}}')^T+o(|\widetilde{\vk{v}}-\widetilde{\vk{v}}'|^2),\quad  \quad \widetilde{\vk{v}},\widetilde{\vk{v}}'\rw\widetilde{\vk{z}}, u\rw\IF.
\EQNY
There exists a non-singular matrix $Q$ such that $Q\Lambda Q^T$ is diagonal, and
set the diagonal is $(c_1,\cdots,c_{n-1})$. Then
\BQNY
\sigma_1(t,Q\widetilde{\vk{v}})=
d-du^{-2}f(t)-\frac{q-2}{2d}\sum_{i=1}^{n-1}c_i(v_i-z_i)^2
+o(|\widetilde{\vk{v}}-\widetilde{\vk{z}}|^2), \widetilde{\vk{v}}\rw\widetilde{\vk{z}},\ u\rw\IF,
\EQNY
and
\BQNY
r_1(t,Q\widetilde{\vk{v}},s,Q\widetilde{\vk{v}}')=1- u^{-2}a(t-s)^\alpha
-\frac{1}{2d}\sum_{i=1}^{n-1}c_i(v_i-z_i)^2+o(|\widetilde{\vk{v}}-\widetilde{\vk{z}}|^2), \widetilde{\vk{v}},\widetilde{\vk{v}}'\rw\widetilde{\vk{z}}, u\rw\IF.
\EQNY
Then set $Y_2(u^{-2/\alpha}t+t_0,\widetilde{\vk{v}})=Y_1(u^{-2/\alpha}t+t_0,Q\widetilde{\vk{v}})$,
defined on a set $[-S_1,S_2]\times(Q^{-1}\widetilde{\mathcal{S}}^\delta_q(1))$. We know that the point $Q\widetilde{\vk{z}}$ is a interior point of $Q^{-1}\widetilde{\mathcal{S}}^\delta_q(1)$.
Then the proof follows by similar arguments as in the proof of \cite{Pit96} [Theorem 8.2]. Consequently, we get
\BQNY
\pk{\sup_{(t,\vk{v})\in \mathcal{D}^1_\delta}\frac{Y(u^{-2/\alpha}t+t_0 ,\vk{v})}{1+u^{-2}f(t)}>u}
&=&\pk{\sup_{(t,\widetilde{\vk{v}})\in [-S_1,S_2]\times(Q^{-1}\widetilde{\mathcal{S}}^\delta_q(1))}Y_2(u^{-2/\alpha}t+t_0, \widetilde{\vk{v}})>u}\\
&\sim&\mathcal{P}_{\alpha,ad^{-2}}^{\frac{1}{d^{2}}f(t)}
[-S_1,S_2]\LT(\prod_{i=1}^{n-1}\mathcal{P}_{2,1}^{(q-2)t^2}
(-\IF,\IF)\RT)\Psi\LT(\frac{u}{d}\RT)\\
&=&\mathcal{P}_{\alpha,ad^{-2}}^{\frac{1}{d^{2}}f(t)}
[-S_1,S_2](2-p)^{(1-n)/2}\Psi\LT(\frac{u}{d}\RT),\ u\rw\IF,
\EQNY
 where we use the fact in \cite{MR1206354} that
 \BQNY
 \mathcal{P}_{2,1}^{(q-2)t^2}(-\IF,\IF)=\sqrt{1+\frac{1}{q-2}}=(2-p)^{-1/2},
 \EQNY
 and
 \BQNY
 \pk{\sup_{(t,\vk{v})\in [-S_1,S_2]\times\mathcal{S}_q}\frac{Y(u^{-2/\alpha}t+t_0 ,\vk{v})}{1+u^{-2}f(t)}>u}\sim 2^n\mathcal{P}_{\alpha,ad^{-2}}^{\frac{1}{d^{2}}f(t)}
[-S_1,S_2](2-p)^{(1-n)/2}\Psi\LT(\frac{u}{d}\RT),\ u\rw\IF.
 \EQNY
 Case 2) $p\in(2,\IF]$ and $M=2m$. Again we need to find the asymptotics of single term in \eqref{sum1}, to wish namely for a maximum point $(0, \vk{v}^1_{+})$, $\vk{v}^1_{+}=(1,0,\cdots,0)$ of variance $\sigma_1^2(t,\vk{v})$.
 hence the fields $\frac{Y(u^{-2/\alpha}t+t_0 ,\vk{v})}{1+u^{-2}f(t)}$ can be represented as
\BQNY
Y_1(u^{-2/\alpha}t+t_0, \widetilde{\vk{v}})
=\sum_{i=2}^{n}v_id_i \frac{X_i(u^{-2/\alpha}t+t_0)}{1+u^{-2}f(t)}
+\LT(1-\sum_{i=2}^{n}\abs{v_i}^q\RT)^{1/q}d_1\frac{X_1(u^{-2/\alpha}t+t_0)}{1+u^{-2}f(t)},
 \widetilde{\vk{v}}=(v_2,\cdots,v_{n}),
\EQNY
which is defined in $[-S_1,S_2]\times\widetilde{\mathcal{S}}^\delta_q(1)$ where
$$\widetilde{\mathcal{S}}^\delta_q(1)
=\LT\{\widetilde{\vk{v}}:\LT(\LT(1-\sum_{i=2}^{n}\abs{v_i}^q \RT)^{1/q},v_2,\cdots,v_{n},\RT)\in\mathcal{S}^\delta_q(1)\RT\},$$
is a small neighborhood of $\widetilde{\vk{0}}:=(0,\cdots,0)\in\R^{n-1}$.
On $[-S_1,S_2]\times\widetilde{\mathcal{S}}^\delta_q(1)$, the variance
\BQNY
\sigma_1^2(t,\widetilde{\vk{v}}):=\frac{1}{(1+u^{-2}f(t))^2}\sigma_1^2(\widetilde{\vk{v}})
:=\frac{1}{(1+u^{-2}f(t))^2}\LT[\sum_{i=2}^{n}d_i^2v_i^2+d_n^2\LT(1-\sum_{i=2}^{n}\abs{v_i}^q\RT)^{2/q}\RT]
\EQNY
of $Y_1(u^{-2/\alpha}t+t_0, \widetilde{\vk{v}})$ attains its maximum $1$ at $(0,\widetilde{\vk{0}})$ where $\widetilde{\vk{0}}$ is a interior point of a set $\widetilde{\mathcal{S}}^\delta_q(1)$. We can write the following Taylor expansion for $\sigma_1(t,\widetilde{\vk{v}})$
\BQNY
\sigma_1(t,\widetilde{\vk{v}})=1-u^{-2}f(t)
-\frac{1}{q}\sum_{i=2}^n\abs{v_i}^q+o\LT(\sum_{i=2}^n|v_i|^q\RT), \widetilde{\vk{v}}\rw\widetilde{\vk{0}},\ u\rw\IF,
\EQNY
and the following expansion for the correlation function $r_1(t,\widetilde{\vk{v}},s,\widetilde{\vk{v}}')$ of $Y_1(u^{-2/\alpha}t+t_0, \widetilde{\vk{v}})$
\BQNY
r_1(t,\widetilde{\vk{v}},s,\widetilde{\vk{v}}')=1- u^{-2}a(t-s)^\alpha
-\frac{1}{2}\sum_{i=2}^nd_i^2(v_i-v_i')^2+o\LT(\sum_{i=2}^nd_i^2(v_i-v_i')^2\RT), \widetilde{\vk{v}},\widetilde{\vk{v}}'\rw\widetilde{\vk{0}}, u\rw\IF.
\EQNY
Then the proof again follows by similar arguments as in the proof of \cite{Pit96} [Theorem 8.2]. Consequently, we get
\BQNY
\pk{\sup_{(t,\vk{v})\in \mathcal{D}^1_\delta}\frac{Y(u^{-2/\alpha}t+t_0 ,\vk{v})}{1+u^{-2}f(t)}>u}
&=&\pk{\sup_{(t,\widetilde{\vk{v}})\in [-S_1,S_2]\times(\widetilde{\mathcal{S}}^\delta_q(1))}Y_1(u^{-2/\alpha}t+t_0, \widetilde{\vk{v}})>u}\\
&\sim&\mathcal{P}_{\alpha,a}^{f(t)}
[-S_1,S_2]\Psi\LT(u\RT),
\EQNY
 and
 \BQNY
 \pk{\sup_{(t,\vk{v})\in [-S_1,S_2]\times\mathcal{S}_q}\frac{Y(u^{-2/\alpha}t+t_0 ,\vk{v})}{1+u^{-2}f(t)}>u}\sim 2m\mathcal{P}_{\alpha,a}^{f(t)}
[-S_1,S_2]\Psi\LT(u\RT), \ u\rw\IF.
 \EQNY
 Case 3) $p=2$.
 By \nelem{dd}, we know that $ \sigma_1^2(t,\vk{v})$ attains its maximum (equal to 1) over $[-S_1,S_2]\times\mathcal{S}_q$ only at points on $\{(0,\vk{v}),\vk{v}\in\mathcal{S}_q,v_i=0,m+1\leq i\leq n\}$.
 The fields $\frac{Y(u^{-2/\alpha}t+t_0 ,\vk{v})}{1+u^{-2}f(t)}$ again can be represented as
\BQNY
Y_1(u^{-2/\alpha}t+t_0, \widetilde{\vk{v}})
=\sum_{i=2}^{n}v_id_i \frac{X_i(u^{-2/\alpha}t+t_0)}{1+u^{-2}f(t)}
+\LT(1-\sum_{i=2}^{n}v_i^q\RT)^{1/q}d_1\frac{X_1(u^{-2/\alpha}t+t_0)}{1+u^{-2}f(t)},
 \widetilde{\vk{v}}=(v_2,\cdots,v_{n}),
\EQNY
which is defined in $[-S_1,S_2]\times\widetilde{\mathcal{S}}_q$ where
$$\widetilde{\mathcal{S}}_q
=\LT\{\widetilde{\vk{v}}:\LT(\LT(1-\sum_{i=2}^{n}v_i^q \RT)^{1/q},v_2,\cdots,v_{n},\RT)\in\mathcal{S}_q\RT\}.$$
On $[-S_1,S_2]\times\widetilde{\mathcal{S}}_q$, the variance
\BQNY
\sigma_1^2(t,\widetilde{\vk{v}}):=\frac{1}{(1+u^{-2}f(t))^2}\sigma_1^2(\widetilde{\vk{v}})
:=\frac{1}{(1+u^{-2}f(t))^2}\LT[\sum_{i=2}^{n}d_i^2v_i^2+d_n^2\LT(1-\sum_{i=2}^{n}v_i^q\RT)^{2/q}\RT]
\EQNY
of $Y_1(u^{-2/\alpha}t+t_0, \widetilde{\vk{v}})$ attains its maximum $1$ at $\LT\{(0,\widetilde{\vk{v}}),\widetilde{\vk{v}}\in\widetilde{\mathcal{S}}_q,v_i=0,m+1\leq i\leq n\RT\}$.  Furthermore, following the arguments as in \cite{Pitchi1994} we conclude that
 $\sigma_1(t,\widetilde{\vk{v}})$  and the correlation function $r_1(t,\vk{v},s,\widetilde{\vk{v}})$ of $Y_1(u^{-2/\alpha}t+t_0, \widetilde{\vk{v}})$ have the following asymptotic expansions:
\BQNY
\sigma_1(t,\widetilde{\vk{v}})=1-u^{-2}f(t)
-\sum_{i=m+1}^n\frac{1-d_i^2}{2}\abs{v_i}^2+o\LT(\sum_{i=m+1}^n\frac{1-d_i^2}{2}\abs{v_i}^2+u^{-2}\RT), \widetilde{\vk{v}}\rw\widetilde{\vk{0}},\ u\rw\IF,
\EQNY
and the following expansion for the correlation function $r_1(t,\widetilde{\vk{v}},s,\widetilde{\vk{v}}')$ of $Y_1(u^{-2/\alpha}t+t_0, \widetilde{\vk{v}})$
\BQNY
r_1(t,\widetilde{\vk{v}},s,\widetilde{\vk{v}}')=1- u^{-2}a(t-s)^\alpha
-\frac{1}{2}\sum_{i=2}^nd_i^2(v_i-v_i')^2+o\LT(\sum_{i=2}^nd_i^2(v_i-v_i')^2+u^{-2}\RT), \widetilde{\vk{v}},\widetilde{\vk{v}}'\rw\widetilde{\vk{0}}, u\rw\IF.
\EQNY
Then the proof follows by similar arguments as in the proof of \cite{PL2015} [Theorem 6.1] with the case $\mu=\nu$. Consequently, we get
\BQNY
\pk{\sup_{(t,\vk{v})\in [-S_1,S_2]\times\mathcal{S}_q}\frac{Y(u^{-2/\alpha}t+t_0 ,\vk{v})}{1+u^{-2}f(t)}>u}
&=&\pk{\sup_{(t,\widetilde{\vk{v}})\in [-S_1,S_2]\times(\widetilde{\mathcal{S}}_q)}Y_1(u^{-2/\alpha}t+t_0, \widetilde{\vk{v}})>u}\\
&\sim&\mathcal{P}_{\alpha,a}^{f(t)}
[-S_1,S_2]\frac{\sqrt{2\pi}2^{\frac{(2-m)}{2}}u^{m-3}}
{\Gamma(m/2)}\LT(\prod_{i=m+1}^n(1-d_i^2)^{-\frac{1}{2}}\RT)\Psi\LT(u\RT).
\EQNY
{\bf Step 2:} Next we proceed to the proof of \eqref{P2p2}.
Setting $a_{u,k}=(a(ku^{-2/\alpha}S+t_0))^{1/\alpha}$, then for any $k\in K_u$ with
$\lim_{u\rw\IF}\sup_{k\in K_u}\abs{k u^{-2/\alpha}}\leq \theta$ and $t\in[-S_1,S_2]$ when $u$ large enough
\BQNY
(a-\vn_\theta)^{1/\alpha} \leq a_{u,k}\leq (a+\vn_\theta)^{1/\alpha}
\EQNY
holds for some $\vn_\theta\in (0,a)$.\\
Then we have
 \BQNY
 \pk{\sup_{t\in [-S_1,S_2]}Z(u^{-\frac{2}{\alpha}}(t+kS)+t_0 ) >u_k}&=& \pk{\sup_{t\in [-a_{u,k}S_1,a_{u,k}S_2]}Z(u^{-\frac{2}{\alpha}}(a_{u,k}^{-1}t+kS)+t_0 )>u_k}\\
 &\leq & \pk{\sup_{t\in [-(a+\vn_\theta)^{1/\alpha}S_1,(a+\vn_\theta)^{1/\alpha}S_2]}Z(u^{-\frac{2}{\alpha}}(a_{u,k}^{-1}t+kS)+t_0 )>u_k}\\
 &=&:\Pi^{+}(u)
 \EQNY
 and
 \BQNY
 \pk{\sup_{t\in [-S_1,S_2]}Z(u^{-\frac{2}{\alpha}}(t+kS)+t_0 ) >u_k}
 &\geq & \pk{\sup_{t\in [-(a-\vn_\theta)^{1/\alpha}S_1,(a-\vn_\theta)^{1/\alpha}S_2]}Z(u^{-\frac{2}{\alpha}}(a_{u,k}^{-1}t+kS)+t_0 )>u_k}\\
 &=&:\Pi^{-}(u).
 \EQNY
 We notice that by assumption (iv)
 \BQNY
 Cov(X(u^{-\frac{2}{\alpha}}(a_{u,k}^{-1}t+kS)+t_0),X(u^{-\frac{2}{\alpha}}kS+t_0))
 &\sim& 1-a(u^{-\frac{2}{\alpha}}kS+t_0)\abs{u^{-\frac{2}{\alpha}}a_{u,k}^{-1}t}^\alpha\\
 &=&1-u^{-2}\abs{t}^{\alpha},\ u\rw\IF.
 \EQNY
For $\Pi^{+}(u)$ and $\Pi^{-}(u)$, when $p=1$, \eqref{P2p2} follows with the same arguments as in {\bf Step 1}.\\
 When $p\in(1,\IF]$, for $\Pi^{+}(u)$ and $\Pi^{-}(u)$ we use the similar arguments as in in {\bf Step 1} with
 $Y_1(u^{-2/\alpha}(a^{-1}_{u,k}t+kS)+t_0 ,\widetilde{\vk{v}})=Y(u^{-2/\alpha}(a^{-1}_{u,k}t+kS)+t_0 ,\widetilde{\vk{v}})$.\\
 When $p\in(1,2)$,
 \BQNY
\sigma_1(t,Q\widetilde{\vk{v}})=
d-\frac{q-2}{2d}\sum_{i=1}^{n-1}c_i(v_i-z_i)^2
+o\LT(|\widetilde{\vk{v}}-\widetilde{\vk{z}}|^2\RT), \widetilde{\vk{v}}\rw\widetilde{\vk{z}},
\EQNY
and
\BQNY
r_1(t,Q\widetilde{\vk{v}},s,Q\widetilde{\vk{v}}')=1- u^{-2}(t-s)^\alpha
-\frac{1}{2d}\sum_{i=1}^{n-1}c_i(v_i-z_i)^2+o\LT(|\widetilde{\vk{v}}-\widetilde{\vk{z}}|^2+u^{-2}\RT), \widetilde{\vk{v}},\widetilde{\vk{v}}'\rw\widetilde{\vk{z}}, u\rw\IF.
\EQNY
When $p\in(2,\IF]$,
\BQNY
\sigma_1(t,\widetilde{\vk{v}})=1
-\frac{1}{q}\sum_{i=2}^n\abs{v_i}^q+o\LT(\sum_{i=2}^n|v_i|^q\RT), \widetilde{\vk{v}}\rw\widetilde{\vk{0}},
\EQNY
and
\BQNY
r_1(t,\widetilde{\vk{v}},s,\widetilde{\vk{v}}')=1- u^{-2}(t-s)^\alpha
-\frac{1}{2}\sum_{i=2}^nd_i^2(v_i-v_i')^2+o\LT(\sum_{i=2}^nd_i^2(v_i-v_i')^2+u^{-2}\RT), \widetilde{\vk{v}},\widetilde{\vk{v}}'\rw\widetilde{\vk{0}}, u\rw\IF.
\EQNY
When $p=2$,
\BQNY
\sigma_1(t,\widetilde{\vk{v}})=1
-\sum_{i=m+1}^n\frac{1-d_i^2}{2}\abs{v_i}^2+o\LT(\sum_{i=m+1}^n\frac{1-d_i^2}{2}\abs{v_i}^2\RT), \widetilde{\vk{v}}\rw\widetilde{\vk{0}},
\EQNY
and
\BQNY
r_1(t,\widetilde{\vk{v}},s,\widetilde{\vk{v}}')=1- u^{-2}(t-s)^\alpha
-\frac{1}{2}\sum_{i=2}^nd_i^2(v_i-v_i')^2+o\LT(\sum_{i=2}^nd_i^2(v_i-v_i')^2+u^{-2}\RT), \widetilde{\vk{v}},\widetilde{\vk{v}}'\rw\widetilde{\vk{0}}, u\rw\IF.
\EQNY
We get that as $u\rw\IF$
\BQNY
&&\Pi^{+}(u)\sim \mathcal{H}_{\alpha}[-S_1(a+\vn_\theta)^{1/\alpha}d^{-2/\alpha},S_2(a+\vn_\theta)^{1/\alpha}d^{-2/\alpha}]
\pk{Z(t_0)) >u_k},\\
&&\Pi^{-}(u)\sim \mathcal{H}_{\alpha}[-S_1(a-\vn_\theta)^{1/\alpha}d^{-2/\alpha},S_2(a-\vn_\theta)^{1/\alpha}d^{-2/\alpha}]
\pk{Z(t_0))>u_k}.
\EQNY
Thus \eqref{P2p2} follows.\\
Further, if letting $\theta\rw 0$ in \eqref{P2p2}, we get \eqref{P2p}.

\QED

\BEL\label{in1}
Assume that Gaussian vector process  $\vk{X}(t)$ with independent marginals which have unit variances, correlation functions $r(t)$ is the same as in \nelem{im1}. Further, set $K_u$ a family of index sets and $u_k$ satisfying that
\BQNY
\lim_{u\rw\IF}\sup_{k\in K_u}\LT|\frac{u_k}{u}-1\RT|=0.
\EQNY
Let $\vn_0$ be such that for all $s,t\in[t_0-\vn_0,t_0+\vn_0]$,
\BQNY
\frac{a}{2}\abs{t-s}^{\alpha}\leq1-r(s,t)\leq2a\abs{t-s}^\alpha.
\EQNY
Then we can find a constant $\mathbb{C}$ such that for all $S>0$ and $ T_2-T_1>S$,
\BQNY
{\limsup}_{u\rw\IF}\sup_{k\in K_u}\frac{\pk{\mathcal{A}_1(u_k),\mathcal{A}_2(u_k)}}
{\pk{Z(t_0)>u_k}}\leq \mathbb{C}\exp\LT(-\frac{a}{8}|T_2-T_1-S|^\alpha\RT),
\EQNY
where $\mathcal{A}_i(u_k)=\{\sup_{t\in [T_i,T_i+S]}Z(u^{-2/\alpha}(t+kS)+t_0)>u_k\},\ i=1,2$, and
\BQNY
\lim_{u\rw\IF}\sup_{k\in K_u}\abs{u^{-2/\alpha}k S}\leq \vn_0.
\EQNY
\EEL
\prooflem{in1}
Through this proof, $\mathbb{C}_i, i\in\N$ are some positive constant.\\
When $p=1$, set $\mathcal{W}=\{\vk{w}=(w_1,\cdots,w_n): w_i=\pm 1,\ i=1,\cdots, n\}$.
We have by \cite{Uniform2016}[Theorem 3.1] for $u$ large enough
\BQNY
&&\pk{\mathcal{A}_1(u_k),\mathcal{A}_2(u_k)}\\
&&=\pk{\sup_{t\in [T_1,T_1+S]}\sum_{i=1}^n\abs{d_i X_i(u^{-2/\alpha}(t+kS)+t_0)}>u_k,
\sup_{s\in [T_2,T_2+S]}\sum_{i=1}^n\abs{d_i X_i(u^{-2/\alpha}(s+kS)+t_0)}>u_k}\\
&&\leq{\sum_{\vk{w}\in\mathcal{W}}}
\pk{\sup_{t\in [T_1,T_1+S]}\sum_{i=1}^nw_id_iX_i(u^{-2/\alpha}(t+kS)+t_0)>u_k,\sup_{s\in [T_2,T_2+S]}\sum_{i=1}^nw_id_iX_i(u^{-2/\alpha}(s+kS)+t_0)>u_k}\\
&&=2^{n}
\pk{\sup_{t\in [T_1,T_1+S]}\sum_{i=1}^nd_iX_i(u^{-2/\alpha}(t+kS)+t_0)>u_k,\sup_{s\in [T_2,T_2+S]}\sum_{i=1}^nd_iX_i(u^{-2/\alpha}(s+kS)+t_0)>u_k}\\
&&=2^{n}
\pk{\sup_{t\in [T_1,T_1+S]}\LT(\sum_{i=1}^nd^2_i\RT)^{1/2}X_1(u^{-2/\alpha}(t+kS)+t_0)>u_k,\sup_{s\in [T_2,T_2+S]}\LT(\sum_{i=1}^nd^2_i\RT)^{1/2}X_1(u^{-2/\alpha}(s+kS)+t_0)>u_k}\\
&&\leq \mathbb{C}_0\exp\LT(-\frac{a}{8}|T_2-T_1-S|^\alpha\RT)\pk{Z(t_0)>u_k}.
\EQNY
When $p\in(1,\IF]$, set $Y_u(t,\vk{v})=\sum_{i=1}^n d_iv_iX_{i}(u^{-2/\alpha^*}(t+kS)+t_0),(t,\vk{v})\in\R\times\mathcal{S}_q$ which is a centered Gaussian field and $\mathcal{S}_q^{\delta}=\{\vk{v}\in\mathcal{S}_q:d^2-\sum_{i=1}^nd_i^2v_i^2\leq \delta\},\delta>0$.\\
Below for $\Delta_1,\Delta_2\subseteq \R^{n+1}$, denote
\BQNY
\mathcal{Y}_u(\Delta_1,\Delta_2)=\pk{\sup_{(t,\vk{v})\in\Delta_1}Y_u(t,\vk{v})>u_k
,\sup_{(t,\vk{v})\in\Delta_2}Y_u(t,\vk{v})>u_k}.
\EQNY
We have
\BQNY
\pk{\mathcal{A}_1(u_k),\mathcal{A}_2(u_k)}&\geq&
\mathcal{Y}_u([T_1,T_1+S]\times\mathcal{S}^\delta_q,[T_2,T_2+S]\times\mathcal{S}^\delta_q),\\
\pk{\mathcal{A}_1(u_k),\mathcal{A}_2(u_k)}&\leq&
\mathcal{Y}_u([T_1,T_1+S]\times\mathcal{S}^\delta_q,[T_2,T_2+S]\times\mathcal{S}^\delta_q)+
\mathcal{Y}_u([T_1,T_1+S]\times\mathcal{S}^\delta_q,[T_2,T_2+S]\times(\mathcal{S}_q
\setminus\mathcal{S}^\delta_q))\\
&&+\mathcal{Y}_u([T_1,T_1+S]\times(\mathcal{S}_q\setminus\mathcal{S}^\delta_q)
,[T_2,T_2+S]\times\mathcal{S}^\delta_q),
\EQNY
and
\BQNY
\mathcal{Y}_u([T_1,T_1+S]\times\mathcal{S}^\delta_q,[T_2,T_2+S]\times(\mathcal{S}_q
\setminus\mathcal{S}^\delta_q))
&\leq& \pk{
\sup_{(t,\vk{v})\in[T_2,T_2+S]\times(\mathcal{S}_q\setminus\mathcal{S}^\delta_q)}Y_u(t,\vk{v})>u_k}\\
&\leq& \exp\LT(-\frac{\LT(u_k-\mathbb{C}_1\RT)^2}{2(d^2-\delta)}\RT)\\
&=&o\LT(\pk{Z(t_0)>u_k}\RT),
\EQNY
as $u\rw\IF$ where the last second inequality follows from Borell inequality and the fact that
\BQNY
\sup_{(t,\vk{v})\in[T_2,T_2+S]\times(\mathcal{S}_q\setminus\mathcal{S}^\delta_q)}
\Var(Y_u(t,\vk{v}))=\sup_{\vk{v}\in(\mathcal{S}_q\setminus\mathcal{S}^\delta_q)}
\LT(\sum_{i=1}^n d_i^2v_i^2\RT)\leq d^2-\delta.
\EQNY
Similarly, we have
\BQNY
\mathcal{Y}_u([T_1,T_1+S]\times(\mathcal{S}_q\setminus\mathcal{S}^\delta_q)
,[T_2,T_2+S]\times\mathcal{S}^\delta_q)=o\LT(\pk{Z(t_0)>u_k}\RT),\ u\rw\IF.
\EQNY
Then we just need to focus on
$$\Pi(u):=\mathcal{Y}_u([T_1,T_1+S]\times\mathcal{S}^\delta_q,[T_2,T_2+S]\times\mathcal{S}^\delta_q).$$
We split $\mathcal{S}^\delta_q$ into sets of small diameters $\{\partial \mathcal{S}_i, 0\leq i\leq \mathcal{N}^*\}$, where
\BQNY
\mathcal{N}^*=\sharp\{\partial\mathcal{S}_i\}<\IF.
\EQNY
Further, we see that $\Pi(u)\leq \Pi_1(u)+\Pi_2(u)$ with
\BQNY
\Pi_1(u)=\underset{\partial\mathcal{S}_i\cap\partial\mathcal{S}_l=\emptyset}
{\sum_{0\leq i,l\leq \mathcal{N}^*}}
\mathcal{Y}_u([T_1,T_1+S]\times\partial\mathcal{S}_i,[T_2,T_2+S]\times\partial\mathcal{S}_l),\\
\Pi_2(u)=\underset{\partial\mathcal{S}_i\cap\partial\mathcal{S}_l\neq\emptyset}{\sum_{0\leq i,l\leq \mathcal{N}^*}}
\mathcal{Y}_u([T_1,T_1+S]\times\partial\mathcal{S}_i,[T_2,T_2+S]\times\partial\mathcal{S}_l),
\EQNY
where ${\partial\mathcal{S}_i\cap\partial\mathcal{S}_l\neq\emptyset}$ means $\partial\mathcal{S}_i,\partial\mathcal{S}_l$ are identical or adjacent, and ${\partial\mathcal{S}_i\cap\partial\mathcal{S}_l=\emptyset}$ means $\partial\mathcal{S}_i,\partial\mathcal{S}_l$ are neither identical nor adjacent.
Denote the distance of two set $\vk{A}, \vk{B}\in\R^n$ as
\BQNY
\rho(\vk{A},\vk{B})=\inf_{\vk{x}\in \vk{A},\vk{y}\in \vk{B}}\|\vk{x}-\vk{y}\|_2.
\EQNY
if $\partial\mathcal{S}_i\cap\partial\mathcal{S}_l=\emptyset$, then there exists some small positive constant $\rho_0$ (independent of $i,l$) such that $\rho(\partial\mathcal{S}_i,\partial\mathcal{S}_l)>\rho_0$.
Next we estimate $\Pi_1(u)$. For any $u\geq 0$
\BQNY
\Pi_1(u)\leq \pk{\underset{\vk{v}\in\partial\mathcal{S}_i,\vk{w}\in\partial\mathcal{S}_i}
{\sup_{(t,s)\in[T_1,T_1+S]\times[T_2,T_2+S]}}Z_u(t,\vk{v},s,\vk{w})>2u_k},
\EQNY
where $Z_u(t,\vk{v},s,\vk{w})=Y_u(t,\vk{v})+Y_u(s,\vk{w}), \ t,s\geq 0, \vk{v},\vk{w}\in\R^n$.\\
When  $u$ is sufficiently large for $(t,s)\in[T_1,T_1+S]\times[T_2,T_2+S],\vk{v}\in\partial\mathcal{S}_i\subset[-2,2]^n,
\vk{w}\in\partial\mathcal{S}_i\subset[-2,2]^n$, with $\rho(\partial\mathcal{S}_i,\partial\mathcal{S}_l)>\rho_0$ we have
\BQNY
Var(Z_u(t,\vk{v},s,\vk{w}))&\leq&\sum_{i=1}^n(v_i^2+w_i^2+2v_iw_i)d_i^2\\
&\leq&4d^2-2\sum_{i=1}^n(v_i-w_i)^2d_i^2\\
&=&4d^2-2d_n^2\rho_0\\
&\leq&d^2(4-\delta_0),
\EQNY
for some $\delta_0>0$. Therefore, it follows from the Borell inequality that
\BQNY
\Pi_1(u)\leq \mathbb{C}_2 \mathcal{N}^*\exp\LT(-\frac{\LT(2u_k-\mathbb{C}_3\RT)^2}{2d^2(4-\delta_0)}\RT)
&=&o\LT(\pk{Z(t_0)>u_k}\RT),
\ u\rw\IF,
\EQNY
with
$$\mathbb{C}_3=\E{\underset{(\vk{v},\vk{w})\in[-2,2]^2n}{\sup_{(t,s)\in[T_1,T_1+S]\times[T_2,T_2+S]}}Z_u(t,\vk{v},s,\vk{w})}<\IF.$$
Now we consider $\Pi_2(u)$. Similar to the argumentation as in {\bf Step1} of the proof of \nelem{im1}. we set $\widetilde{Y}_u(t,\widetilde{\vk{v}})=Y_u(t,Q\widetilde{\vk{v}})$ and $\widetilde{Z}_u(t,\widetilde{\vk{v}},s,\widetilde{\vk{w}})=\widetilde{Y}_u(t,\widetilde{\vk{v}})+\widetilde{Y}_u(s,\widetilde{\vk{w}})$ with $\widetilde{\vk{v}},\widetilde{\vk{w}}\in \R^{n-1}$.
Since for $(t,s)\in[T_1,T_1+S]\times[T_2,T_2+S],\widetilde{\vk{v}}\in[-2,2]^{n-1},
\widetilde{\vk{w}}\in[-2,2]^{n-1}$, we have
\BQNY
2d^2\leq Var(\widetilde{Z}_u(t,\widetilde{\vk{v}},s,\widetilde{\vk{w}}))&\leq&\sum_{i=1}^n(v_i^2+w_i^2+2r(u^{-2/\alpha}(t+kS)+t_0,
u^{-2/\alpha}(s+kS)+t_0)v_i w_i)d_i^2\\
&\leq&2d^2+2\LT(1-\frac{a}{2}u^{-2}\abs{t-s}^\alpha\RT)\sum_{i=1}^n v_iw_id_i^2\\
&\leq&4d^2-d^2au^{-2}\abs{t-s}^\alpha\\
&\leq&4d^2-d^2au^{-2}\abs{T_2-T_1-S}^\alpha.
\EQNY
Set
\BQNY
\overline{Z}_u(t,\widetilde{\vk{v}},s,\widetilde{\vk{w}})=
\frac{\widetilde{Z}_u(t,\widetilde{\vk{v}},s,\widetilde{\vk{w}})}
{Var(Z_u(t,\widetilde{\vk{v}},s,\widetilde{\vk{w}}))}.
\EQNY
Borrowing the arguments of the proof in \cite{Pit96} [Lemma 6.3] we show that
\BQNY
\E{\LT(\overline{Z}_u(t,\widetilde{\vk{v}},s,\widetilde{\vk{w}})-\overline{Z}_u(t',\widetilde{\vk{v}'},s',\widetilde{\vk{w}'})\RT)}
\leq 4\LT(\E{(\widetilde{Y}_u(t,\widetilde{\vk{v}})-\widetilde{Y}_u(t',\widetilde{\vk{v}'}))^2}
+\E{(Y_u(s,\widetilde{\vk{w}})-Y_u(s',\widetilde{\vk{w}'}))^2}\RT).
\EQNY
Moreover,  since
when $p\in(1,2)$,
\BQNY
r_1(t,Q\widetilde{\vk{v}},s,Q\widetilde{\vk{v}}')=1- u^{-2}a(t-s)^\alpha
-\frac{1}{2d}\sum_{i=1}^{n-1}c_i(v_i-z_i)^2+o\LT(|\widetilde{\vk{v}}-\widetilde{\vk{z}}|^2+u^{-2}\RT), \widetilde{\vk{v}},\widetilde{\vk{v}}'\rw\widetilde{\vk{z}}, u\rw\IF.
\EQNY
When $p\in(2,\IF)$,
\BQNY
r_1(t,\widetilde{\vk{v}},s,\widetilde{\vk{v}}')=1- u^{-2}a(t-s)^\alpha
-\frac{1}{2}\sum_{i=2}^nd_i^2(v_i-v_i')^2+o\LT(\sum_{i=2}^nd_i^2(v_i-v_i')^2+u^{-2}\RT), \widetilde{\vk{v}},\widetilde{\vk{v}}'\rw\widetilde{\vk{0}}, u\rw\IF.
\EQNY
When $p=2$,
\BQNY
r_1(t,\widetilde{\vk{v}},s,\widetilde{\vk{v}}')=1- u^{-2}a(t-s)^\alpha
-\frac{1}{2}\sum_{i=2}^nd_i^2(v_i-v_i')^2+o\LT(\sum_{i=2}^nd_i^2(v_i-v_i')^2+u^{-2}\RT), \widetilde{\vk{v}},\widetilde{\vk{v}}'\rw\widetilde{\vk{0}}, u\rw\IF.
\EQNY

Then we have
\BQNY
\E{(Y_u(t,\widetilde{\vk{v}})-Y_u(t',\widetilde{\vk{v}'}))^2}\leq 4d^2a u^{-2}\abs{t-t'}^\alpha+2\sum_{i=2}^n(v_i-v_i')^2.
\EQNY
Therefore
\BQNY
\E{\LT(\overline{Z}_u(t,\widetilde{\vk{v}},s,\widetilde{\vk{w}})
-\overline{Z}_u(t',\widetilde{\vk{v}'},s',\widetilde{\vk{w}'})\RT)}\leq 16d^2a u^{-2}\abs{t-t'}^\alpha+16d^2a u^{-2}\abs{s-s'}^\alpha+8\sum_{i=2}^n(v_i-v_i')^2
+8\sum_{i=2}^n(w_i-w_i')^2.
\EQNY
Set $\zeta(t,s,\widetilde{\vk{v}},\widetilde{\vk{w}}), t,s\geq0, \vk{v},\vk{w}\in\R^{n-1}$ is a stationary Gaussian field with unit variance and correlation function
\BQNY
r_{\zeta}(t,s,\widetilde{\vk{v}},\widetilde{\vk{w}})=\exp\LT(-9 d^2at^\alpha-9 d^2as^\alpha-5\sum_{i=2}^nv_i^2-5\sum_{i=2}^nw_i^2\RT).
\EQNY
Then
\BQNY
\Pi_2(u)&\leq& \pk{\underset{\widetilde{\vk{v}}\in Q^{-1}\mathcal{S}_q,\widetilde{\vk{w}}\in Q^{-1}\mathcal{S}_q}
{\sup_{(t,s)\in[T_1,T_1+S]\times[T_2,T_2+S]}}\widetilde{Z}_u(t,\widetilde{\vk{v}},s,\widetilde{\vk{w}})>2u_k}\\
&\leq&\pk{\underset{\widetilde{\vk{v}}\in Q^{-1}\mathcal{S}_q,\widetilde{\vk{w}}\in Q^{-1}\mathcal{S}_q}
{\sup_{(t,s)\in[T_1,T_1+S]\times[T_2,T_2+S]}}\zeta(u^{-2/\alpha}t,u^{-2/\alpha}s,\widetilde{\vk{v}},\widetilde{\vk{w}})
>\frac{2u_k}{\sqrt{4d^2-d^2au^{-2}\abs{T_2-T_1-S}^\alpha}}}.
\EQNY
Then following the similar argumentation as in \cite{EnkelejdJi2014Chi}, we have
\BQNY
\Pi_2(u)\leq \mathbb{C}_4u_k^{M-2}\exp\LT(-\frac{u_k^{2}}{2d^2}-\frac{a}{8}\abs{T_2-T_1-S}^\alpha\RT)
\EQNY
where $M=0$ when $p\in(1,2)\cup(2,\IF]$ and $M=m$ when $p=2$.
Thus we have
\BQNY
\limsup_{u\rw\IF}\frac{\Pi_2(u)}{\pk{Z(t_0)>u_k}}\leq \mathbb{C}_5\exp\LT(-\frac{a}{8}|T_2-T_1-S|^\alpha\RT).
\EQNY
Thus we complete the proof.
\COM{
Set $Y(t,\vk{v})=\sum_{i=1}^n d_iv_iX_{i}(t),(t,\vk{v})\in\R\times\mathcal{S}_q$ which is a centered Gaussian field.
First we have for $u$ large
\BQNY
&&2\leq4-4a|t-s|^\alpha\leq\Var(X_i(s)+X_i(t))=2+2r(|t-s|)\leq4-a|t-s|^\alpha
\leq4-a\abs{T_2-T_1-S}\alpha u^{-2/c},\\
&&1-Cov(X_i(s)+X_i(t),X_i(s')+X_i(t'))\leq 1-r(|s-s'|)+1-r(|t-t'|)\leq 2a|s-s'|^{\alpha}+2a|t-t'|^{\alpha},
\EQNY
with $s\in[t_0+T_1u^{-2/\alpha^*},t_0+(T_1+S)u^{-2/\alpha^*}],
t\in[t_0+T_2u^{-2/\alpha^*},t_0+(T_2+S)u^{-2/\alpha^*}]$.\\
Set $\mathcal{S}_q=\{\vk{v}\in\R^n:||\vk{v}||_q=1\}, \mathcal{S}'_q=\{\vk{v}'\in\R^n:||\vk{v}'||_q=1\}$. We divide $\mathcal{S}_q$ and $\mathcal{S}_q$ respectively into $2^n$ parts according to $v_i\geq 0$ and $v_i\leq 0$  which are denote by $\mathcal{S}_q(j),\mathcal{S}'_q(j), j=1,\cdots,2^n$.
We know
\BQNY
\pk{\mathcal{A}_1(u_k),\mathcal{A}_2(u_k)}
&=&\pk{\sup_{(t,\vk{v})\in[T_1,T_1+S]\times\mathcal{S}_q}Y(u^{-2/\alpha^*}t+t_0,\vk{v})>u^{1/c}_k
,\sup_{(s,\vk{v}')\in[T_2,T_2+S]\times\mathcal{S}'_q}Y(u^{-2/\alpha^*}s+t_0,\vk{v}')>u^{1/c}_k}\\
&\leq&\pk{\sup_{(t,s,\vk{v},\vk{v}')\in[T_1,T_1+S]\times[T_2,T_2+S]\times\mathcal{S}_q
\times\mathcal{S}'_q}Y(u^{-2/\alpha^*}t+t_0,\vk{v})+Y(u^{-2/\alpha^*}s+t_0,\vk{v}')>2u^{1/c}_k}\\
&\leq&\sum_{j=1}^{2^n}\sum_{k=1}^{2^n}\pk{\sup_{(t,s,\vk{v},\vk{v}')\in[T_1,T_1+S]\times[T_2,T_2+S]\times\mathcal{S}_q(j)
\times\mathcal{S}'_q(k)}Y(u^{-2/\alpha^*}t+t_0,\vk{v})+Y(u^{-2/\alpha^*}s+t_0,\vk{v}')>2u^{1/c}_k}.
\EQNY
Next we consider one term in the sum of the upper inequality.
Without lose of generality, we set $\mathcal{S}_q(1)=\{\vk{v}\in\R^n:||\vk{v}||_q=1, v_i\geq 0\}$ and $\mathcal{S}'_q(1)=\{\vk{v}'\in\R^n:||\vk{v}'||_q=1, v'_i\geq 0\}$.}

\QED

\proofExamp{EX1} Note that the variance function $\sigma^2(t)$ of  $B_\alpha(t)$ attain its maximum over $[0,1]$ at $t=1$ and
\BQN\label{ex1aa}
\sigma(t)\sim 1-\frac{\alpha}{2}(1-t),\ \ r(s,t)\sim 1-\frac{1}{2}\abs{s-t}^{\alpha}, \ s,t \uparrow 1.
\EQN
For $g(t)=-(1-t)^{1/2}, \ t\in[0,1]$, by \netheo{Thm2} with $c=1$ we get the result.

\proofprop{riskmodel} Note that the variance function $\sigma^2(t)$ of $B_\alpha(t)$ attains its maximum over $[0,1]$ at $t=1$ and
\eqref{ex1aa} is satisfied.
Since
\BQNY
\pk{\inf_{t\in[0,1]}U(t)<0}=\pk{\sup_{t\in[0,1]}\LT(\sum_{i=1}^n \abs{d_iB_\alpha^{i}(t)}^2+w(1-t)\RT)>u+w},
\EQNY
for $g(t)=w(1-t), \ t\in[0,1]$, by \netheo{Thm2} and {\bf Remarks} \ref{rem11} ii) with $c=2$ and $p=2$ we get the result.

\QED

\proofprop{OUCSP}
Note that $V_i(t),\ 1\leq i\leq n$ are stationary with unit variance and correlation function $r(s,t)$ satisfies
\BQNY
 r(s,t)\sim 1-2\abs{s-t}, \abs{s-t}\rw 0, \quad \text{and} \quad r(s,t)<1, \forall s\neq t.
\EQNY
By \netheo{Thm0} with $c=2$ and $p=2$ we get the result.

\QED

{\bf Acknowledgement}:Thanks  to  the  referees  for  their  comments  and  suggestions  which  significantly improved the manuscript. Thanks to Prof. Enkelejd Hashorva for his suggestions. Thanks to  Swiss National Science Foundation Grant no.  200021-166274.
\bibliographystyle{plain}
\bibliography{LP}
\end{document}